\documentclass[12pt]{amsart}

\usepackage[a4paper]{geometry}
\usepackage{latexsym,enumerate}

\usepackage{amsmath,amsthm,amsopn,amstext,amscd,amsfonts,amssymb}

\usepackage{stackrel}
\usepackage[thicklines]{cancel}
\allowdisplaybreaks
\usepackage{color}

\newcommand{\R}{\mathbb{R}}
\newcommand{\N}{\mathbb{N}}

\newcommand{\car}{{\raise0pt\hbox{{\LARGE $\chi$}}}}
\newcommand{\sg}{{\rm \; sign}}

\newcommand{\Div}{\hbox{\rm div\,}}
\newcommand{\dis }{{\mathcal D}' }
\newcommand{\z }{{\bf z}}

\newenvironment{pf}{\noindent{\sc Proof}.\enspace}{\rule{2mm}{2mm}\medskip}
\newtheorem{Theorem}{Theorem}[section]
\newtheorem{Corollary}[Theorem]{Corollary}
\newtheorem{Definition}[Theorem]{Definition}
\newtheorem{Lemma}[Theorem]{Lemma}
\newtheorem{Proposition}[Theorem]{Proposition}

\newtheorem{Remark}[Theorem]{Remark}

\newcommand{\norma}[2]{\|#1\|_{\lower 4pt \hbox{$\scriptstyle #2$}}}
\newcommand{\ldoble}{\langle\!\langle}
\newcommand{\rdoble}{\rangle\!\rangle}

.

\def\estrella{\buildrel\ast\over\rightharpoonup }
\newcommand{\h}{{\mathcal H}^{N-1}}

\begin{document}

\title[The inhomogeneous 1-Laplace evolution equation]{Existence and uniqueness for the inhomogeneous 1-Laplace evolution equation revisited}

\author[M. Latorre and  S. Segura de Le\'on]
{Marta Latorre and  Sergio Segura de Le\'on}

\address{Marta Latorre: Matem\'atica Aplicada, Ciencia e Ingenier\'ia de los Materiales y Tecnolog\'ia Electr\'onica, Universidad Rey Juan Carlos
\hfill\break\indent C/Tulip\'an s/n 28933, M\'ostoles, Spain
\hfill\break\indent ORCID iD: 0000-0001-9859-3809}\email{\tt marta.latorre@urjc.es}

\address{Sergio Segura de Le\'on: Departament d'An\`{a}lisi Matem\`atica,
Universitat de Val\`encia,
\hfill\break\indent Dr. Moliner 50, 46100 Burjassot, Spain
\hfill\break\indent ORCID iD: 0000-0002-8515-7108}
\email  {\tt sergio.segura@uv.es }

\thanks{}


\bigskip
\begin{abstract}
In this paper we deal with an inhomogeneous  parabolic Dirichlet problem involving the 1-Laplacian operator. We show the existence of a unique solution when data belong to $L^1(0,T;L^2(\Omega))$ for every $T>0$. As a consequence, global existence and uniqueness for data in $L^1_{loc}(0,+\infty;L^2(\Omega))$ is obtained. Our analysis retrieves previous results in a correct and complete way.
\end{abstract}

\maketitle

\section{Introduction}
The aim of this paper is to prove the existence of a unique solution to the following evolution problem:
\begin{equation}\label{prob}
\left\{\begin{array}{rcll}
 \displaystyle u'  -  \Delta_1 u & \!=\!& f(t,x)\,&\hbox{in }(0,T)\times\Omega\,, \\[2mm]
 u &\!=\! &0\,  & \hbox{on } (0,T)\times\partial\Omega \,,\\[2mm]
 u(0,x)&\!=\!&u_{0}(x)& \hbox{in }\Omega\,,
\end{array}\right.
\end{equation}
where $\Omega$ is a bounded open set in $\R^N$ ($N\ge2$) with Lipschitz boundary and $T>0$.
Henceforth, the sign $'$ stands for the derivative with respect to time variable $t$ while  $\Delta_1=\Div\left(\frac{Du}{|Du|}\right)$ is the so--called 1--Laplacian operator.
As far as data are concerned, we will take an initial datum $u_0\in L^2(\Omega)$ and a source $f\in L^1(0,T;L^2(\Omega))$.

The homogeneous problem, $f\equiv 0$, was solved in \cite{ABCM} (see also \cite{ACDM, ACM}). Using nonlinear semigroup theory, authors were able to introduce a concept of solution and to prove existence and uniqueness. Since the natural space to analyze the stationary problem is the space $BV(\Omega)$ of functions of bounded variation, one of their crucial tasks was to make clear the quotient $\frac{Du}{|Du|}$, for $u\in BV(\Omega)$, even if $Du$ vanishes. The successful method was to consider a bounded vector field $\z\in L^\infty(\Omega; \R^N)$ which plays the role of that quotient in the sense that it satisfies $\|\z\|_\infty\le 1$ and $(\z , Du)=|Du|$ as measures. We stress that the definition of $(\z , Du)$ relies on the Anzellotti pairing theory. Actually, $(\z , Du)$ is a Radon measure which becomes the dot product $\z\cdot\nabla u$ when $u\in W^{1,1}(\Omega)$. The inhomogeneous case was addressed using the techniques of nonlinear semigroup theory in \cite{T} for data $f\in L^2((0,T)\times \Omega)$.

In recent years, a new approach to problem \eqref{prob} has been developed. Following \cite{BDM1}, it was applied to parabolic problems involving the 1--Laplacian in \cite{BDM2} (see also \cite{BDS, BDSS}). Its main feature is that it uses a purely variational approach to deal with time dependent problems, yielding the existence of global parabolic minimizers. The Anzellotti pairings are not used and only the total variation operator is required. A comparison between the two methods and a proof that both approaches lead to the same solutions (under natural assumptions) can be found in \cite{KS}.

In addition to \cite{T}, the inhomogeneous problem is studied in \cite{SW} in a different way. The proof does not lie on nonlinear semigroup theory, although Anzellotti's theory is also involved. Instead of applying the Crandall-Liggett generation theorem, the solution is obtained by approximation in two stages. Firstly one gets a solution for every source belonging to $L^2((0,T)\times\Omega)$ (the same setting of \cite{T}) and then goes beyond studying sources in $L^1(0,T;L^2(\Omega))$. Nevertheless, we consider that some steps of the proof need further explanation. Let us briefly explain why.

In \cite[Proposition 2.9]{SW} is stated that existence of a time derivative $\xi$ of $u$ satisfying $\xi=\Div \z+f$ and Green's formula
imply that $t\mapsto\int_\Omega u(t)^2\,dx$ defines an absolutely continuous function.
In its proof it is assumed that if $\eta\in C_0^\infty(0,T)$ is nonnegative, then
\begin{equation}\label{lim0}
\lim_{\varepsilon\to0}\int_0^T\frac1\varepsilon \int_{t-\varepsilon}^t\left(\eta(s)\int_\Omega (\z(t, x), Du(s,x))\right)\, ds\, dt
\end{equation}
 exists and is equal to
 \[\int_0^T\eta(t)\int_\Omega (\z(t,x), Du(t,x))\, dt\,.\]
 This fact is valid when $u\in L^1(0,T;  W^{1,1}(\Omega))$, as shown in Remark \ref{RIntro}, but it is not justified for a general Radon measure $(\z(t), Du(s))$. The reason is that this device is not a real dot product and we cannot split it into factors.
Moreover, both this product and the total variation $|Du(s)|$ are just Radon measures and not $L^1$--functions.
Thus, the proof is flawed and we are not able to correct it. Instead, we must polish the argument to guarantee similar results.

Our aim in the present paper is to improve \cite{SW}, correcting mistakes and providing more details in order to regain the same results. We take as a starting point solutions to \eqref{prob} having sources in $L^2((0,T)\times\Omega))$. This is so because, in this setting, solutions satisfy $u'(t)\in L^2(\Omega)$ for almost all $t\in(0,T)$ and consequently \cite[Proposition 2.9]{SW} does hold. Hence, we begin by writing this proposition in a way suitable to be extended and so handle more general data (see Proposition \ref{Prop.2.9} below). Next, we use solutions having $L^2$--sources as approximate solutions and translate their main features to any limit solution.

Regarding the proof of existence, we follow the structure of \cite[Theorem 5.1]{SW}, but being more accurate. Several arguments are also taken from \cite{SW} (actually some of them go back to \cite{ABCM}) since we have to deal with similar difficulties. Indeed, it is easy to check the existence of the vector field $\z\in L^\infty((0,T)\times\Omega;\R^N)$ but it is not obvious that $\Div \z$ is in a suitable space to apply Anzellotti's theory since we are just able to see that $\Div \z$ belongs to $L^1(0,T; BV(\Omega)\cap L^2(\Omega))^*$. Following \cite{ABCM}, we have to rewrite this theory in our setting to define the Anzellotti pairings and the trace on the boundary of the normal component of $\z$ as well as to see that a Green's formula holds. In addition, we also need that the limit of approximate solutions can be taken as a test function in problem \eqref{prob} in order to show that it is, in fact, a solution. Furthermore, the identity $(\z(t), Du(t))=|Du(t)|$ which holds as measures for almost all $t\in (0,T)$ is not longer as easy to check and we must analyze a limit similar to \eqref{lim0} studying the Lebesgue points of a related function (see Remark \ref{ident} below).

A remark on uniqueness is also in order. The fact that function given by $t\mapsto\int_\Omega u(t)^2\,dx$ is absolutely continuous for every solution $u$ is essential in the uniqueness proof. Now we cannot infer this property from the identity  $\xi=\Div \z+f$ satisfied by the time derivative of the solution. Hence, to get uniqueness we have to consider condition \eqref{cond2} in our concept of solution.

Summarizing, the proof of \cite[Proposition 2.9]{SW} is not correct but the main results of \cite{SW} hold true with slight modifications in the statements. Nevertheless, more cumbersome arguments are needed to get them.

\bigskip

The plan of this paper is the following. Section 2 is devoted to introduce the notation used through this paper and some preliminary results. In Section 3, we introduce the suitable notion of solution to problem \eqref{prob} and state our starting point. In Section 4  we prove the main result on existence of this paper while in the next section we show the uniqueness of solution. Finally, in an Appendix, we have collected some results from real analysis that are essential over the paper.

\section{Preliminaries}

We use this Section to introduce the notation and some preliminary results that will be used over the course of this paper.

\subsection{Notation}

Henceforth, $T$ will always be a fixed positive number.
We denote by $\Omega$ a bounded open subset in $\R^N$, with $N\ge2$. We also requiere that $\Omega$ has a Lipschitz boundary and $\nu$ will be the outer unit normal vector on $\partial\Omega$ a.e. in $\h(\partial\Omega)$, where $\h$ stands for the $(N-1)$-dimensional Hausdorff measure.

As usual, $L^q(\Omega)$ and $W^{1,q}(\Omega)$ express the Lebesgue and Sobolev spaces respectively (see, for instance, \cite{Br} or \cite{Ev}). Given a Banach space $X$, the symbol  $L^{q}(0,T; X)$  denotes the space of vector valued functions which are strongly measurable and $q$--summable. For instance, we say that $v\in L^q(0,T;L^p(\Omega))$ if $v\colon(0,T)\times\Omega\to\R$ is Lebesgue measurable and
\begin{equation*}
\int_0^T\left(\int_\Omega |v(t,x)|^p\,dx\right)^\frac{q}{p}\,dt<+\infty\,.
\end{equation*}
To simplify the notation, in this case we often write $v(t)$ instead of $v(t,x)$.

\noindent
Throughout this paper, notation $' = \dfrac d{dt}$ will be
used with the meaning $u'=u_t$, the derivative of $u$ with respect to $t$. For this derivative we mean both the derivative in the sense of distributions and one of its extensions (see Definition \ref{T.Der} below).
On the other hand, the symbol $\Div$ denotes the divergence taken in the spatial variables.

\noindent
We also make use of the truncation function defined as follows:
\begin{equation*}
T_k(s)=\min\{|s|,k\}\sg(s)\quad\text{ for all } \;\; s\in\R\,.
\end{equation*}

\subsection{Functions of bounded variation}

In what follows, $BV(\Omega)$ will be the set of all integrable functions in $\Omega$ whose distributional gradient is a Radon measure with finite total variation. For  $v\in BV(\Omega)$ the total variation of its gradient will be denoted as $\int_\Omega|Dv|$.
In general, we denote by $\displaystyle \int_\Omega \varphi |Dv|$ the integral of $\varphi$ with respect to the measure $|Dv|$.
We recall that $BV(\Omega)$ is a Banach space with the norm
\begin{equation*}
\|v\|=\int_\Omega |Dv| +\int_{\Omega} |v|\,dx\,.
\end{equation*}

On the other hand, the notion of a trace on the boundary can be
extended to functions $v \in BV (\Omega)$, so that we may write $v\big|_{\partial\Omega}$, by means of
a surjective bounded operator $BV (\Omega)\to L^1(\partial\Omega)$. As a consequence, an equivalent norm on $BV (\Omega)$ can be defined:
\begin{equation*}
\|v\|_{BV(\Omega)}=\int_\Omega |Dv| +\int_{\partial\Omega} |v|\,d \mathcal H^{N-1}\,.
\end{equation*}

When taking limits, we also make use of the lower semicontinuity of the total variation with respect to the convergence in $L^1(\Omega)$.

We refer to \cite{AFP} for further information on BV-functions.

\noindent
In the present paper, we will widely use the space $BV(\Omega)\cap L^2(\Omega)$, which is a Banach space with the norm defined by
\begin{equation*}
\|v\|=\max\{\|v\|_{BV(\Omega)}, \| v\|_{L^2(\Omega)}\}\,.
\end{equation*}
The dual pairing of $BV(\Omega)\cap L^2(\Omega)$ and its dual will be denoted as
\begin{equation*}
\langle\, \zeta, v\,\rangle_\Omega\quad\text{ with}\quad \zeta\in (BV(\Omega)\cap L^2(\Omega))^*,\;\; v\in BV(\Omega)\cap L^2(\Omega)\,.
\end{equation*}
In the case $\xi=\zeta+f$, for $\zeta\in (BV(\Omega)\cap L^2(\Omega))^*$ and $f\in L^2(\Omega)$, we will write
\begin{equation*}
\ldoble\, \xi, v\,\rdoble_\Omega=\langle \,\zeta, v\,\rangle_\Omega+\int_\Omega f v\, dx\qquad\forall v\in BV(\Omega)\cap L^2(\Omega)\,.
\end{equation*}

\subsection{Anzellotti's theory}

Following \cite{ABCM}, we will
use a vector field $\z\in L^\infty(\Omega;\R^N)$ satisfying $\|\z\|_\infty\le 1$ to play the role of the quotient $\frac{Du}{|Du|}$, even if $|Du|$ vanishes in a zone of positive measure.
If $\z\in L^\infty(\Omega;\R^N)$ with $\Div \z\in L^2(\Omega)$ and $v\in BV(\Omega)\cap L^2(\Omega)$, Anzellotti (see \cite{An}) defined the pairing
\begin{equation*}
\langle\,(\z,Dv),\omega\rangle= - \int_\Omega v\, \omega\, \Div \z\, dx - \int_\Omega v\,\z\cdot \nabla \omega\,dx
\end{equation*}
for all $\omega \in C_0^\infty(\Omega)$. Under these conditions, the pairing $(\z,Dv)$ is actually a Radon measure and it also holds
\begin{equation*}
|(\z,Dv)|\le \|\z\|_\infty |Dv| \qquad\text{ as measures in } \Omega\,.
\end{equation*}
We will write the integration of $\varphi$ with respect to this measure as $\displaystyle \int_\Omega \varphi (\z, Dv)$.

Anzellotti also defined a weak trace on the boundary $\partial\Omega$ of the normal component of the vector field $\z$ (denoted by $[\z,\nu]\in L^\infty(\partial\Omega)$) and proved that the following Gauss-Green formula holds
\begin{equation*}
\int_\Omega v\,\Div \z \,dx + \int_\Omega (\z,Dv) = \int_{\partial \Omega} v [\z,\nu]\,d\h\quad\text{ for all }\;v\in BV(\Omega)\cap L^2(\Omega)\,.
\end{equation*}

\begin{Remark}\label{acop}\rm
For solutions to problem \eqref{prob} having sources in $L^2((0,T)\times\Omega)$ the associated vector field satisfies $\Div\z(t)\in L^2(\Omega)$ for almost all $t\in(0,T)$, so that the Anzellotti theory applies. Nevertheless, for a general source $f\in L^1(0,T; L^2(\Omega))$, we cannot expect that solutions satisfy this property. Hence, we will consider both the pairing $(\z(t), Dv)$ for every $v\in BV(\Omega)\cap L^2(\Omega)$ and the weak trace on the boundary $[\z(t), \nu]$  throughout the existence proof and then check that Green's formula holds. Hereafter, the symbol $(\z,Dv)$, where $\Div \z\in (BV(\Omega)\cap L^2(\Omega))^*$ and $v\in BV(\Omega)\cap L^2(\Omega)$, stands for the distribution defined as
\begin{equation}\label{dis1}
 \langle (\z,Dv), \omega\rangle = -\langle\,\Div\z,v\omega\,\rangle_\Omega - \int_\Omega v\z\cdot \nabla\omega\,dx \,,
\end{equation}
 for all $\omega\in C_0^\infty(\Omega)$.
 \end{Remark}

\subsection{Integration of vector--valued functions}

In this subsection we collect the theory of vector integration needed to our purposes. For more information, we refer to \cite{DU}.

Let $X$ be a Banach space. A function $f\colon(0,T)\to X$ is said to be weakly measurable if for each $x^*\in X^*$ the real function $\langle x^*,f\rangle\colon(0,T)\to \R$ is Lebesgue measurable.

We denote by $L_{w}^{1}(0,T; BV(\Omega))$ the space of all weakly measurable maps
$$
v\colon[0,T]\longrightarrow BV(\Omega)
$$
such that
$\int_{0}^{T}\|v(t)\|_{BV(\Omega)}\, dt<\infty$. Obviously, $L^{1}(0,T; BV(\Omega))\subset L_{w}^{1}(0,T; BV(\Omega))$.

If $f\colon[0,T]\to X$ is weakly measurable such that $\langle x^*,f\rangle\in L^1(0,T)$ for all $x^*\in X^*$, then $f$ is called Dunford integrable. The Dunford integral of $f$ over a measurable set $E\subset (0,T)$ is written
$\int_E f(t)\, dt\in X^{**}$ and given by
\[\left(\int_E f(t)\, dt\right)(x^*)=\int_E \langle x^*,f(t)\rangle\, dt\,.\]

In the case $\int_E f(t)\, dt\in X$ for every measurable set $E\subset (0,T)$, $f$ is called Pettis integrable.

Obviously, both concepts coincide when $X$ is reflexive. Since $BV(\Omega)$ is not reflexive, these concepts may be different in our framework.

In \cite[Lemmas 3-4]{ABCM}, given $u\in L_w^1(0,T; BV(\Omega)\cap L^2(\Omega))$ and $\eta\in C_0^\infty(0,T)$, the function $s\mapsto \eta(s)u(s)$ is proven to be weakly measurable. Then, for every $t\in(0,T)$ and every $0<\varepsilon<t$, the integral
\begin{equation}\label{def_psi}
  \Psi_\varepsilon(t)=\frac{1}{\varepsilon}\int_{t-\varepsilon}^t\eta(s)u(s)\,ds
\end{equation}
is defined as a Dunford integral; it satisfies the following features:
\begin{enumerate}
\item it is actually Pettis integrable,
\item $\Psi_\varepsilon \in C([0,T];BV(\Omega))$,
\item $\Psi_\varepsilon (t)\in BV(\Omega)\cap L^2(\Omega)$ for every $t\in [0,T]$.
\end{enumerate}
Furthermore, if $u\in L^\infty(0,T;L^2(\Omega))$, it is easy to see that
$\|\Psi_\varepsilon\|_{L^\infty(0,T;L^2(\Omega))}\le \|\eta\|_\infty\|u\|_{L^\infty(0,T;L^2(\Omega))}$
and so
\begin{equation}\label{DP-int}
\Psi_\varepsilon\in L^\infty(0,T;BV(\Omega)\cap L^2(\Omega))\,.
\end{equation}

\subsection{The time derivative}
In our equation, both $u'$ and $\Div \z$ refer to derivatives in the sense of distributions.
Since we need to take test functions related with the solution, extensions of these derivatives will be used. This subsection is devoted to clarifying what our definition of time derivative is.

\begin{Definition}\label{defi 1}
Let  $\Psi \in L^1(0,T; BV(\Omega))\cap L^{\infty}(0,T;
L^2(\Omega))$.

We say that $\Psi$ admits a weak derivative in $L_{w}^{1}(0,T;
BV(\Omega))\cap L^{1}(0,T; L^2(\Omega))$ if there exists
$\Theta\in L_{w}^{1}(0,T; BV(\Omega))\cap L^{1}(0,T; L^2(\Omega))$ such
that $\Psi(t) =\int_{0}^{t}\Theta (s) ds$, where the integral is
taken as a Pettis integral.
\end{Definition}

Having in mind the previous subsection, we have that the function $\Psi_\varepsilon$ defined in \eqref{def_psi}
admits a weak derivative in $L_{w}^{1}(0,T; BV(\Omega))\cap L^{1}(0,T; L^2(\Omega))$. Its weak derivative is given by
$\Theta(t)=\dfrac{1}{\varepsilon}\left[\eta(t)u(t)-\eta(t-\varepsilon)u(t-\varepsilon)\right]$.

\begin{Definition}\label{T.Der}
Let $u\in C([0,T];L^2(\Omega))\cap L_w^1(0,T;BV(\Omega))$.

We say that  $\xi \in L^1(0,T; BV(\Omega)\cap L^2(\Omega))^{*}+L^1(0,T;L^2(\Omega))$ is the time derivative of $u$ if
\begin{equation}\label{T-D}
\int_{0}^{T}\ldoble \xi(t),\Psi(t) \rdoble_\Omega dt=-\int_{0}^{T}\int_{\Omega}u(t)\Theta (t)\, dx\, dt\,,
\end{equation}
for every function with compact support in time
$$\Psi \in L^1(0,T; BV(\Omega))\cap L^{\infty}(0,T;L^2(\Omega))$$
 which admits a weak derivative $\Theta\in L_{w}^{1}(0,T; BV(\Omega))\cap L^{1}(0,T; L^2(\Omega))$.
\end{Definition}

\begin{Remark}\rm
We explicitly remark that the dual $L^1(0,T; BV(\Omega)\cap L^2(\Omega))^{*}$ can be identified with $L_{w^*}^\infty(0,T; (BV(\Omega)\cap L^2(\Omega))^*)$, the space of all  $\zeta\colon[0,T]\to (BV(\Omega)\cap L^2(\Omega))^*$ which are *-weakly measurable (that is, the function $t\mapsto\langle\zeta(t),v\rangle$ is measurable for every $v\in BV(\Omega)\cap L^2(\Omega)$) and satisfy $\|\zeta(t)\|_{(BV(\Omega)\cap L^2(\Omega))^*}\in L^\infty(0,T)$ (see \cite{Sch}).
Two functions $\zeta_1$ and $\zeta_2$ are identified if $\langle \zeta_1(t), v\rangle_\Omega=\langle \zeta_2(t), v\rangle_\Omega$ a.e. for all $v\in BV(\Omega)\cap L^2(\Omega)$.
 The duality is given by
\begin{equation*}
\langle\zeta,v\rangle=\int_0^T\langle\zeta(t),v(t)\rangle_\Omega\, dt\,,
\end{equation*}
for all $\zeta\in L_{w^*}^\infty(0,T; (BV(\Omega)\cap L^2(\Omega))^*)$ and all $v\in L^1(0,T;BV(\Omega)\cap L^2(\Omega))$. Consequently the left hand side in \eqref{T-D} makes sense.
\end{Remark}

It is worth specifying how this time derivative works. Consider the function $\Psi_\varepsilon$ defined in \eqref{def_psi}. Then
$\Psi_\varepsilon \in L^\infty(0,T; BV(\Omega) \cap L^2(\Omega))$ and
\begin{equation*}
\int_0^T\ldoble \xi(t), \Psi_\varepsilon(t)\rdoble_\Omega\, dt =- \int_0^T\int_\Omega u(t,x)\big[\eta(t)u(t,x)-\eta(t-\varepsilon)u(t-\varepsilon, x)\big]\, dx\, dt\,.
\end{equation*}

\section{Our starting point}
In this Section, we introduce the definition we will use of solution to our problem and state auxiliary results that hold when the source belongs to $L^2((0,T)\times\Omega)$.

\begin{Definition}\label{def0}
Let $f\in L^1(0,T;L^2(\Omega))$ and $u_0\in L^2(\Omega)$. We say that $u$ is a solution to problem
\begin{equation*}\label{P}
\left\{\begin{array}{rcll}
 \displaystyle u'  -  \Delta_1 u & =& f(t,x)\,&\hbox{in }(0,T)\times\Omega\,, \\[2mm]
 u &= &0\,  & \hbox{on } (0,T)\times\partial\Omega \,,\\[2mm]
 u(0,x)&=&u_{0}(x)& \hbox{in }\Omega\,,
\end{array}\right.
\end{equation*}
if $u\in L_w^1(0,T;BV(\Omega))\cap C([0,T];L^2(\Omega))$ and there exist
\begin{enumerate}
\item[(i)]  $\z\in L^\infty((0,T)\times\Omega)$ such that $\|\z\|_\infty\le 1$ and its distributional divergence $\Div\z$ can be extended to $L^1(0,T;BV(\Omega)\cap L^2(\Omega))^*$ in such a way that, for almost all $t\in(0,T)$, satisfies
\begin{itemize}
\item the weak trace of the normal component of $\z(t)$  is well--defined and verifies $\|[\z(t),\nu]\|_\infty\le1$;
\item $(\z(t), D v)$ is a Radon measure for every $v\in BV(\Omega)\cap L^2(\Omega)$ (recall \eqref{dis1});
\item the following Green's formula is fulfilled:
\[\langle\, \Div \z(t),v\,\rangle_\Omega + \int_\Omega (\z(t),Dv) = \int_{\partial\Omega} v[\z(t),\nu]\,d\h\]
for every $v\in BV(\Omega)\cap L^2(\Omega)$.
\end{itemize}
\item[(ii)] $\xi\in L^1(0,T;BV(\Omega)\cap L^2(\Omega))^*+L^1(0,T;L^2(\Omega))$, which is  the time derivative of $u$ in the sense of Definition \ref{T.Der}.
\end{enumerate}

Moreover,
\begin{equation}\label{cond1}
\xi=\Div \z+f
\end{equation}
and conditions
\begin{align}
&\left(\int_\Omega \frac12 u(t)^2\, dx\right)'+\int_\Omega |Du(t)|+\int_{\partial\Omega}|u(t)|\, d\mathcal H^{N-1}=\int_\Omega f(t)u(t)\, dx\,,\label{cond2}\\[2mm]
&[\z(t),\nu]\in \sg(-u(t)) \quad \h\text{a.e. on }\partial\Omega\,,\label{cond3}
\end{align}
hold for almost every $t\in(0,T)$.
\end{Definition}

\begin{Remark}\label{test_u1}\rm
Applying the above Green's formula pointwise to $u$, we obtain the following identity for almost all $t\in (0,T)$:
 \[\langle \Div\z(t), u(t)\rangle_\Omega=-\int_\Omega (\z(t), Du(t))+ \int_{\partial\Omega}u(t) [\z(t), \nu]\, d\mathcal H^{N-1}\]
 Therefore,
 \begin{equation}\label{cond4}
 \ldoble \xi(t), u(t)\rdoble_\Omega = -\int_\Omega (\z(t), Du(t))+\int_{\partial\Omega}u(t) [\z(t), \nu]\, d\mathcal H^{N-1}+\int_\Omega f(t) u(t) \, dx
 \end{equation}
 holds for almost all $t\in (0,T)$, so that, in some sense, we may take $u$ as test function. Notice, however, that the measurability of these functions are not guaranteed since we do not have $u\in L^1(0,T; BV(\Omega))$.
\end{Remark}

\begin{Remark}\label{test_u0}\rm
It is worth remarking that condition \eqref{cond2} can be written as
\begin{equation*}
\left(\int_\Omega \frac12 u(t)^2\, dx\right)'=
-\int_\Omega |Du(t)|-\int_{\partial\Omega}|u(t)|\, d\mathcal H^{N-1}+\int_\Omega f(t)u(t)\, dx
\end{equation*}
for almost all $t\in(0,T)$. Notice that, owing to $u\in L^1(0,T;BV(\Omega))\cap L^\infty(0,T; L^2(\Omega))$ and $f\in L^1(0,T; L^2(\Omega))$, it follows that the right hand side belongs to $L^1(0,T)$. As a consequence, the function $t\mapsto \int_\Omega  \frac12 u(t)^2\, dx$ is absolutely continuous.
\end{Remark}

\begin{Remark}\label{test_u}\rm
Observe that conditions \eqref{cond2}--\eqref{cond3}, jointly with \eqref{cond4} and Green's formula imply
\begin{multline*}
\left(\int_\Omega \frac12 u(t)^2\, dx\right)'+
\int_\Omega |Du(t)|+\int_{\partial\Omega}|u(t)|\, d\mathcal H^{N-1}
=\int_\Omega f(t) u(t)\, dx\\
=\ldoble \xi(t), u(t)\rdoble_\Omega-
\langle \Div\z(t), u(t)\rangle_\Omega\\
=\ldoble \xi(t), u(t)\rdoble_\Omega+\int_\Omega (\z(t), Du(t))-\int_{\partial\Omega}u(t) [\z(t), \nu]\, d\mathcal H^{N-1}\\
=\ldoble \xi(t), u(t)\rdoble_\Omega+\int_\Omega (\z(t), Du(t))+\int_{\partial\Omega}|u(t)|\, d\mathcal H^{N-1}
\end{multline*}
for almost all $t\in(0,T)$. Hence,
\begin{equation}\label{ec:00}
\left(\int_\Omega \frac12 u(t)^2\, dx\right)'+
\int_\Omega |Du(t)|=\ldoble \xi(t), u(t)\rdoble_\Omega+\int_\Omega (\z(t), Du(t))
\end{equation}
holds for almost all $t\in(0,T)$. This identity suggests that
\[\left(\int_\Omega \frac12 u(t)^2\, dx\right)'=\ldoble \xi(t), u(t)\rdoble_\Omega\]
and
\[\int_\Omega |Du(t)|=\int_\Omega (\z(t), Du(t))\]
but we are not able to check it (see, however, Remark \ref{ident} below).
\end{Remark}

\begin{Proposition}\label{para_unic}
Let  $u$ be a solution to problem \eqref{prob} with time derivative $\xi$. For every nonnegative $\eta\in C_0^\infty(0,T)$,
it is verified that
\begin{equation}\label{unic:0}
\lim_{\varepsilon\to0}\int_0^T\Big\langle\!\!\Big\langle \xi(t), \frac1{\varepsilon}\int_{t-\varepsilon}^t\eta(s) u(s) \, ds\Big\rangle\!\!\Big\rangle_\Omega\, dt=
-\int_0^T\eta'(t)\int_\Omega  \frac12 u(t)^2\, dx\, dt\,,
\end{equation}
and
\begin{equation}\label{unic:-1}
\lim_{\varepsilon\to0}\int_0^T \frac1{\varepsilon}\int_{t-\varepsilon}^t\eta(s) \int_\Omega (\z(t), Du(s)) \, ds\, dt=
\int_0^T\eta(t)\int_\Omega |Du(t)|\, dt\,.
\end{equation}

\end{Proposition}
\begin{pf} Let $\displaystyle  \Psi_\varepsilon(t)=\frac{1}{\varepsilon}\int_{t-\varepsilon}^t\eta(s)u(s)\,ds$. We stress that the integral
\begin{multline*}
\int_0^T\ldoble\, \xi(t), \Psi_\varepsilon(t)\rdoble_\Omega\, dt \\
= \int_0^T \Big\langle\, \Div \z(t),\frac1{\varepsilon}\int_{t-\varepsilon}^t\eta(s) u(s) \, ds \,\Big\rangle_\Omega\,dt + \int_0^T\int_\Omega f(t) \frac1{\varepsilon}\int_{t-\varepsilon}^t\eta(s) u(s) \, ds\,dx\,dt
\end{multline*}
is well--defined since $\Div \z(t) \in L^1(0,T;BV(\Omega)\cap L^2(\Omega))^*$,  function
\[t\mapsto \Psi_\varepsilon(t)=\frac1{\varepsilon}\int_{t-\varepsilon}^t\eta(s) u(s) \, ds\]
belongs to
$L^\infty(0,T;BV(\Omega)\cap L^2(\Omega))$ (recall \eqref{DP-int} above) and $f\in L^1(0,T; L^2(\Omega))$. Here $\frac1{\varepsilon}\int_{t-\varepsilon}^t\eta(s) u(s) \, ds$ is actually a Pettis integral.

Since $\xi$ is the time derivative of $u$, we obtain
\begin{align*}
\int_0^T\ldoble \xi(t), \Psi_\varepsilon (t) \rdoble_\Omega\, dt&=-\frac1{\varepsilon}
\int_0^T \int_\Omega \big(\eta(t)u(t)-\eta(t-\varepsilon) u(t-\varepsilon)\big)u(t) \, dx\, dt\\
&=-\frac1{\varepsilon}\int_0^T \eta(t)\int_\Omega u(t)^2\, dx \, dt+\frac1{\varepsilon}\int_0^T \eta(t-\varepsilon)\int_\Omega u(t-\varepsilon) u(t)\, dx\, dt\\
&=\frac1{\varepsilon} \int_0^T \eta(t) \int_\Omega u(t)\big(u(t+\varepsilon)-u(t)\big)\, dx\, dt\,.
\end{align*}
Observe that inequality $\displaystyle u(t)u(t+\varepsilon)\le \frac12 u(t+\varepsilon)^2+\frac12 u(t)^2$ implies
\begin{equation*}
u(t)\big(u(t+\varepsilon)-u(t)\big)\le \frac12 u(t+\varepsilon)^2-\frac12 u(t)^2\,.
\end{equation*}
Hence
\begin{multline*}
\int_0^T\ldoble \xi(t), \Psi_\varepsilon (t) \rdoble_\Omega\, dt\le \frac1{\varepsilon} \int_0^T\int_\Omega \eta(t)  \Big(\frac12 u(t+\varepsilon)^2-\frac12 u(t)^2\Big)\, dx\, dt\\
=\frac1{\varepsilon} \int_0^T\int_\Omega  \frac12 \Big( \eta(t-\varepsilon)u(t)^2-\eta(t) u(t)^2\Big)\, dx\, dt=\frac1{\varepsilon} \int_0^T\big( \eta(t-\varepsilon)-\eta(t)\big)\int_\Omega  \frac12 u(t)^2\, dx\, dt\,.
\end{multline*}
Next, we will compute $\int_0^T\ldoble \xi(t), \Psi_\varepsilon (t) \rdoble_\Omega\, dt$ in a different way. Indeed, we apply Green's formula to deduce
\begin{multline}\label{ec:13.00}
-\frac1{\varepsilon} \int_0^T\big( \eta(t-\varepsilon)-\eta(t)\big)\int_\Omega  \frac12 u(t)^2\, dx\, dt\le-\int_0^T\ldoble \xi(t), \Psi_\varepsilon (t) \rdoble_\Omega\, dt\\
=-\int_0^T \frac1{\varepsilon}\int_{t-\varepsilon}^t \eta(s)\Big[\langle\, \Div\z(t), u(s) \,\rangle_\Omega +\int_\Omega u(s)  f (t)\, dx\Big]\, ds\, dt\\
=\int_0^T \frac1{\varepsilon}\int_{t-\varepsilon}^t \eta(s)\Big[\int_\Omega  (\z(t), Du(s))
-\int_{\partial\Omega} u(s) [\z(t), \nu]\, d\mathcal H^{N-1} -\int_\Omega u(s)f(t) \, dx\Big]\, ds\, dt\\
\le \int_0^T \frac1{\varepsilon}\int_{t-\varepsilon}^t \eta(s)\Big[\int_\Omega  |Du(s)|
+\int_{\partial\Omega} |u(s)| \, d\mathcal H^{N-1} -\int_\Omega u(s)f(t) \, dx\Big]\, ds\, dt\,.
\end{multline}
Our next concern is to take the limit on the right hand side of \eqref{ec:13.00} as $\varepsilon$ tends to 0. We deal with the first and second terms thanks to Lemma \ref{lema0}; it is enough to consider the functions given as
$\eta(s)\int_\Omega |Du(s)|$ and $\eta(s)\int_{\partial\Omega}|u(s)|\, d\mathcal H^{N-1}$. The remainder is handled as follows
\begin{align*}
\left|\frac1{\varepsilon}\int_{t-\varepsilon}^t \eta(s)\int_\Omega u(s) f(t)  \, dx\, ds\right| & \le  \frac1{\varepsilon}\int_{t-\varepsilon}^t \eta(s)\left(\int_\Omega  u(s)^2\, dx\right)^{\frac12}\left(\int_\Omega  f(t)^2\, dx\right)^{\frac12}\, ds\\
 & \le   \|f(t)\|_{L^2(\Omega)}\frac1{\varepsilon}\int_{t-\varepsilon}^t \eta(s)\left(\int_\Omega  u(s)^2\, dx\right)^{\frac12}ds\\
 & \le  \|\eta\|_\infty  \|u\|_{L^\infty(0,T;L^2(\Omega))}\|f(t)\|_{L^2(\Omega)}
\end{align*}
which belongs to $L^1(0,T)$.
Hence, letting $\varepsilon$ go to 0 in \eqref{ec:13.00}, it yields
\begin{multline}\label{parte}
 \int_0^T \eta'(t)\int_\Omega  \frac12 u(t)^2\, dx\, dt\le-\lim_{\varepsilon\to0}\int_0^T\ldoble \xi(t), \Psi_\varepsilon (t) \rdoble_\Omega\, dt\\
 \le
 \int_0^T \eta(t)\Big[\int_\Omega  |Du(t)|
+\int_{\partial\Omega} |u(t)|\, d\mathcal H^{N-1}-\int_\Omega u(t)f(t) \, dx\Big]\, dt\,.
\end{multline}

On the other hand, $u$ is a solution to problem \eqref{prob}. So, recalling Remark \ref{test_u} and inserting that identity in \eqref{parte}, we obtain
\[
 \int_0^T \eta'(t)\int_\Omega  \frac12 u(t)^2\, dx\, dt\le-\lim_{\varepsilon\to0}\int_0^T\ldoble \xi(t), \Psi_\varepsilon (t) \rdoble_\Omega\, dt
 \le
 \int_0^T \eta'(t)\int_\Omega  \frac12 u(t)^2\, dx\, dt\,,
\]
from where \eqref{unic:0} follows.
Going back to \eqref{ec:13.00} and letting $\varepsilon$ go to 0, inequatilies become equalities and it follows from condition \eqref{cond3} that limit \eqref{unic:-1} holds.
\end{pf}

\begin{Remark}\label{ident}\rm
Taken into account Lemma \ref{lema1}, we may understand \eqref{unic:-1} as follows: value $\eta(t)\int_\Omega |Du(t)|$ is ``almost'' the approximate limit of $\eta(s)\int_\Omega (\z(t), Du(s))$ at $s=t$ for every $\eta\in C_0^\infty(0,T)$. (Note, however, that limit is not pointwise but in mean.) So, in some sense, this result allows us to identify $\int_\Omega (\z(t), Du(t))=\int_\Omega |Du(t)|$ for almost all $t\in (0,T)$.
\end{Remark}

\subsection{Source data in $L^2((0,T)\times\Omega)$ }
In order to prove the existence of solution to our problem, we need a previous result which appears in \cite{SW}, even though we will restrict our analysis to data $f \in L^2((0, T)\times\Omega)$.

By \cite[Theorem 4.1]{SW}, for each $f \in L^2((0, T)\times\Omega)$ and each $u_0 \in L^2(\Omega)$ there exists a solution to \eqref{prob} which satisfies $u\in C([0,T];L^2(\Omega))\cap L^1(0,T;BV(\Omega))$ as well as $u' \in L^2((\delta, T)\times \Omega)$ for all $\delta>0$.
So, $\Div \z(t)\in L^2(\Omega)$ for almost all $t\in(0,T)$ and Anzellotti's theory applies.

As a consequence of the Green formula, it is obtained that $u'$ is the time derivative of $u$ in the sense of Definition \ref{T.Der}.
Indeed, for every $v\in L^1(0,T; BV(\Omega)\cap L^2(\Omega))$, we get
\begin{equation*}
\int_0^T\langle \Div \z (t), v(t)\rangle_\Omega dt=-\int_0^T\int_\Omega(\z(t), Dv(t))\, dt+\int_0^T\int_{\partial\Omega}v(t)\, [\z(t), \nu]\, d\mathcal H^{N-1} dt
\end{equation*}
which is well-defined and so $u'\in L^1(0,T; BV(\Omega)\cap L^2(\Omega))^*+L^2(0,T; L^2(\Omega))$. The condition \eqref{T-D} is now easy to check.

Thus, to see that $u$ is a solution to problem \eqref{prob} in the sense of Definition \ref{def0}, it just remains to show that satisfies $\ldoble \xi(t), u(t)\rdoble_\Omega=\frac12\left( \int_\Omega u(t)^2\,dx\right)^\prime$ for almost all $t\in (0,T)$. We check this condition in the following result.
We will also apply Proposition \ref{Prop.2.9} in the proof of Theorem \ref{exis}.

\begin{Proposition}\label{Prop.2.9}
For any $u\in C([0,T]; L^2(\Omega))$ satisfying $u'\in L^2((\delta,T)\times\Omega)$ for every $\delta>0$, the following identity holds for almost all $t\in (0,T)$:
\begin{equation*}
\dfrac{1}{2}\left(\int_\Omega u(t)^2\,dx\right)' = \int_\Omega u'(t)u(t)\,dx\,.
\end{equation*}
\end{Proposition}
\begin{pf} We will check that the proof of \cite[Proposition 2.9]{SW} works in this case.

Let $\eta\in C_0^\infty(0,T)$ and let $\varepsilon>0$ be small enough to perform the following calculations. Observe that  $u' \eta \in L^2((0, T)\times \Omega)$.

Making the same calculations that at the beginning of \cite[Proposition 2.9]{SW} we get
\begin{multline}\label{ec.prop.01}
-\int_0^T\int_\Omega\dfrac{\eta(t-\varepsilon)-\eta(t)}{-\varepsilon} \dfrac{u(t)^2}{2}\,dx\,dt
 = \dfrac{1}{2}\int_0^T\int_\Omega\dfrac{u(t+\varepsilon)-u(t)}{\varepsilon}u(t+\varepsilon)\eta(t)\,dx\,dt\\
 +\dfrac{1}{2}\int_0^T\int_\Omega\dfrac{u(t+\varepsilon)-u(t)}{\varepsilon}u(t)\eta(t)\,dx\,dt
=\dfrac{1}{2}(I_1+I_2)\,.
\end{multline}
We now consider the auxiliary function $\displaystyle\Psi_\varepsilon(t)=\frac{1}{\varepsilon}\int_{t-\varepsilon}^t\eta(s)u(s)\,ds$, which is a Pettis integral. Since $\Psi_\varepsilon\in L^\infty(0,T; L^2(\Omega))$, it satisfies
$\ldoble\,u'(t),\Psi_\varepsilon(t)\,\rdoble_\Omega\in L^1(0,T)$. In addition, \cite[Proposition 2.9]{SW} also yields $I_2=\int_0^T\ldoble\, u'(t),\Psi_\varepsilon(t)\,\rdoble_\Omega\,dt$.

We also consider $\displaystyle\Phi_\varepsilon(t)=\dfrac{1}{\varepsilon}\int_{t-\varepsilon}^t\eta(s)u(s+\varepsilon)\,ds$, which is a Pettis integral. Then $\ldoble\, u'(t),\Phi_\varepsilon(t)\,\rdoble_\Omega\in L^1(0,T)$ and it can be proved that
 $\displaystyle I_1=\int_0^T\ldoble\, u'(t), \Phi_\varepsilon(t)\,\rdoble_\Omega\,dt$.

 Therefore, from \eqref{ec.prop.01} we deduce that
\begin{multline}\label{ec.prop.03}
-\int_0^T\int_\Omega \dfrac{\eta(t-\varepsilon)-\eta(t)}{-\varepsilon}\dfrac{u(t)^2}{2}\,dx\,dt \\
= \dfrac{1}{2}\left(\int_0^T\ldoble\,u'(t),\Psi_\varepsilon(t)\,\rdoble_\Omega\,dt+\int_0^T\ldoble\,u'(t), \Phi_\varepsilon(t)\,\rdoble_\Omega\,dt\right)\,.
\end{multline}
Now, we take limits when $\varepsilon$ tends to $0$.
 The first term on the right hand side is handled as follows.
\begin{multline*}
|\ldoble\, u'(t),\Psi_\varepsilon(t)\,\rdoble_\Omega | \le \int_\Omega  |\Psi_\varepsilon(t)||u'(t)|\,dx \\
\le \left(\int_\Omega \left(\dfrac{1}{\varepsilon}\int_{t-\varepsilon}^t |\eta(s)||u(s)|\,ds\right)^2\,dx\right)^\frac{1}{2}\left(\int_\Omega \chi_{{supp\,}\eta}|u'(t)|^2\,dx\right)^\frac{1}{2}\\
  \le \left( \dfrac{1}{\varepsilon} \int_{t-\varepsilon}^t \eta(s)^2 \int_\Omega|u(s)|^2\,dx\,ds\right)^\frac{1}{2}\left(\int_\Omega \chi_{{supp\,}\eta}|u'(t)|^2\,dx\right)^\frac{1}{2}\,,
\end{multline*}
and both factors belong to $L^2(0,T)$.
Moreover, applying  Lemma A.1, we get that the family
\begin{equation*}
\dfrac{1}{\varepsilon} \int_{t-\varepsilon}^t \eta(s)^2 \left(\int_\Omega|u(s)|^2\,dx\right)\,ds
\end{equation*}
converges in $L^1(0,T)$, so that the generalized dominated convergence theorem and the pointwise convergence
$\displaystyle\lim_{\varepsilon\to0}\frac1{\varepsilon}\int_{t-\varepsilon}^t \eta(s) u(s, x)\,ds = \eta(t) u(t, x)$ imply
\begin{equation*}
\ldoble\, u'(t),\Psi_\varepsilon(t)\,\rdoble_\Omega\to \int_\Omega \eta(t) u(t) u'(t)\, dx\,,
\end{equation*}
where the convergence holds in $L^1(0,T)$. Hence,
\begin{equation}\label{ec.prop.025}
\lim_{\varepsilon\to0}\int_0^T  \ldoble\, u'(t), \Psi_\varepsilon(t)\,\rdoble_\Omega\, dt=\int_0^T \int_\Omega \eta(t) u(t) u'(t)\, dx\, dt\,.
\end{equation}
Similarly, we deduce that
\begin{equation}\label{ec.prop.026}
\lim_{\varepsilon\to0}\int_0^T
\ldoble\, u'(t), \Phi_\varepsilon(t)\,\rdoble_\Omega\, dt=\int_0^T \int_\Omega \eta(t) u(t) u'(t)\, dx\, dt\,.
\end{equation}

Letting $\varepsilon$ go to $0$ in \eqref{ec.prop.03}, by \eqref{ec.prop.025} and \eqref{ec.prop.026}, we get
\begin{equation*}
- \int_0^T\eta'(t)\int_\Omega \dfrac{u(t)^2}{2}\,dx\,dt = \int_0^T\int_\Omega \eta(t)u(t)u'(t)\,dx\,dt\,.
\end{equation*}
Since this identity holds for every $\eta\in C_0^\infty(0,T)$, it yields
\begin{equation*}
\left(\dfrac{1}{2}\int_\Omega u(t)^2\,dx\right)' = \int_\Omega u(t)u'(t)\,dx\,,
\end{equation*}
for almost all $t\in (0,T)$.
\end{pf}

\begin{Remark}\label{RProp.2.9}\rm
Reasoning as in the proof of Proposition \ref{Prop.2.9}, we can also prove that
\begin{equation*}
\dfrac{1}{2}\left(\int_\Omega u(t)^2\omega\,dx\right)' = \int_\Omega u'(t)u(t)\,\omega\,dx\,
\end{equation*}
 for all $\omega\in C_0^\infty(\Omega)$.
\end{Remark}

Notice that solutions with source in $L^2((0,T)\times\Omega)$ are unique. To see it, just argue as in the uniqueness proof of \cite{SW} and take into account that the function
\begin{equation*}
t\longmapsto \int_\Omega u(t)^2\,dx\quad\text{ is absolutely continuous in } (0,T)\,.
\end{equation*}

\section{Existence of solution if $f \in L^1(0,T;L^2(\Omega))$}
This section is devoted to prove the existence of a solution to problem \eqref{prob}. We will follow the proof of Theorem 5.1 in \cite{SW}, but trying to guarantee all details.

\begin{Theorem}\label{exis} If $f\in L^1(0,T;L^2(\Omega))$ and $u_0\in L^2(\Omega)$, then there exists, at least, a solution to problem \eqref{prob}.
\end{Theorem}
\begin{pf}
Since $f\in L^1(0,T;L^2(\Omega))$, there exists a sequence $\{ f_n\}\subset L^2(0,T;L^2(\Omega))$ such that $f_n\to f$ in $L^1(0,T;L^2(\Omega))$. Furthermore, each approximating problem
\begin{equation}\label{Pn}
\left\{\begin{array}{rcll}
 \displaystyle u_n^\prime  -  \Delta_1 u_n & =& f_n(t,x)\,&\hbox{in }(0,T)\times\Omega\,,\\[2mm]
 u_n &= &0\,  & \hbox{on } (0,T)\times\partial\Omega \,,\\[2mm]
 u_n(x,0)&=&u_{0}(x)& \hbox{in }\Omega\,,
\end{array}\right.
\end{equation}
has a solution $u_n\in L_w^1(0,T;BV(\Omega))\cap C([0,T];L^2(\Omega))$ whose time derivative  satisfies $u_n^\prime\in L^1(0,T;BV(\Omega)\cap L^2(\Omega))^*+L^1(0,T;L^2(\Omega))$ and $u_n^\prime \in L^2((\delta, T)\times\Omega)$ for  all $\delta>0$. Moreover, there exists a vector field $\z_n\in L^\infty((0,T)\times\Omega;\R^N)$ with $\|\z_n\|_\infty\le 1$ such that
\begin{enumerate}
\item $u_n^\prime(t) =\Div\z_n(t)+f_n(t)$ in $\dis(\Omega)$,
\item $(\z_n(t),Du_n(t)) = |Du_n(t)|$ as measures in $\Omega$,
\item $[\z_n(t),\nu]\in \sg(-u_n(t))$,
\item $\int_\Omega u'_n(t)u_n(t)\, dx=\frac12\left( \int_\Omega u_n(t)^2\,dx\right)^\prime$\,,
\end{enumerate}
holds for almost every $t\in (0,T)$. This last identity is due to Proposition \ref{Prop.2.9}.

\medskip
\noindent
Our purpose is to check that the sequence $\{u_n\}$ converges to a function $u$, which is a solution to problem \eqref{prob}.
We divide the proof in several steps.

\subsection{Step 1: A priori estimates}

 We begin applying \eqref{cond2} and  H\"older's inequality to get
\begin{multline*}
 \frac12 \left(\int_\Omega u_n(t)^2\,dx\right)^\prime +\int_\Omega|Du_n(t)|+\int_{\partial \Omega} |u_n(t)|\,d\h = \int_\Omega f_n(t)\,u_n(t)\,dx \\
 \le \left( \int_\Omega f_n(t)^2\,dx\right)^\frac{1}{2}\left( \int_\Omega u_n(t)^2\,dx \right)^\frac{1}{2} \,.
\end{multline*}
Integrating now between $0$ and $t\in(0,T]$, it yields
\begin{multline*}
\dfrac{1}{2}\int_\Omega u_n(t)^2\,dx  - \dfrac{1}{2}\int_\Omega u_n(0)^2\,dx+\int_0^t\int_\Omega|Du_n(s)|\,ds+\int_0^t\int_{\partial \Omega}|u_n(s)|\,d\h \,ds \\
 \le \int_0^t\left( \int_\Omega f_n(s)^2\,dx\right)^\frac{1}{2}\left( \int_\Omega u_n(s)^2\,dx \right)^\frac{1}{2} \,ds\,.
\end{multline*}
Denoting $\sigma_n(t)=\|u_n(t)\|_{L^2(\Omega)}$ and disregarding nonnegative terms, the previous inequality becomes
\begin{equation*}
\sigma_n(t)^2\le \sigma_n(0)^2 + 2 \int_0^t\|f_n(s)\|_{L^2(\Omega)}\, \sigma_n(s)\,ds\,.
\end{equation*}
Now, due to Gronwall's Lemma of \cite{W1} (see also \cite{W2} for a slightly extension) and the fact that all approximating problems have the same initial data, we get
\begin{equation*}
\sigma_n(t)\le\sigma_n(0) + \int_0^t\|f_n(s)\|_{L^2(\Omega)} \le \| u_0\|_{L^2(\Omega)} + \|f_n\|_{L^1(0,T;L^2(\Omega))}
\end{equation*}
which is bounded because $f_n\to f$ in $L^1(0,T;L^2(\Omega))$. So there exists a constant $C_1>0$ such that
\begin{equation}\label{ec:1}
\| u_n(t)\|_{L^2(\Omega)}\le\| u_0\|_{L^2(\Omega)} + \|f_n\|_{L^1(0,T;L^2(\Omega))} \le C_1\quad\text{ for all }\; t \in [0,T]\,.
\end{equation}
Moreover,
\begin{multline*}
\int_\Omega u_n(t)^2\,dx  +2\int_0^t\Big[\int_\Omega|Du_n(s)|+\int_{\partial \Omega}|u_n(s)|\,d\h \Big]\,ds \\
 \le\int_\Omega u_n(0)^2\,dx+ 2\int_0^t \|f_n(s)\|_{L^2(\Omega)}\| u_n(s)\|_{L^2(\Omega)} \,ds\\
 \le\int_\Omega u_0^2\,dx+ 2\,C_1\,\|f_n\|_{L^1(0,T;L^2(\Omega))} = C_2
\end{multline*}
for all $t\in(0,T)$, and so
\begin{equation}\label{mainEstimate}
\max_{t\in[0,T]}\int_\Omega u_n(t)^2\,dx +2\int_0^T\Big[\int_\Omega|Du_n(s)|+\int_{\partial \Omega}|u_n(s)|\,d\h\Big] \,ds \le C_2\,.
\end{equation}

\subsection{Step 2: Convergence of the sequence $(u_n)_n$ in $L^\infty(0,T; L^2(\Omega))$}

We next check that $\{u_n\}_{n=1}^{\infty}$ is a Cauchy sequence. We already know that $\{u_n\}$ is bounded in the space $L^\infty(0,T;L^2(\Omega))$ owing to \eqref{ec:1}.

\noindent
Taking $u_n(t)-u_m(t)$ as a test in problem \eqref{Pn} for $n$ and then taking it for $m$ lead to
\begin{multline*}
\int_{\Omega}u^\prime_n(t)(u_{n}(t)-u_m(t))\, dx +\int_{\Omega}(\z_n(t),D(u_{n}(t)-u_m(t)) ) \\
- \int_{\partial\Omega}(u_{n}(t)-u_m(t))[\z_n(t),\nu]\, d\mathcal H^{N-1} =\int_{\Omega}f_n(t)(u_{n}(t)-u_m(t))\, dx
\end{multline*}
and
\begin{multline*}
\int_{\Omega}u^\prime_m(t)(u_{n}(t)-u_m(t)) \, dx+\int_{\Omega}(\z_m(t),D(u_{n}(t)-u_m(t)) )\\
 - \int_{\partial\Omega}(u_{n}(t)-u_m(t))[\z_m(t),\nu]\, d\mathcal H^{N-1}
 =\int_{\Omega}f_m(t)(u_{n}(t)-u_m(t))\, dx\,.
\end{multline*}
Subtracting both expressions yields
\begin{multline*}
\int_{\Omega}(u_{n}(t)-u_m(t))\,(u_{n}(t)-u_m(t))^\prime\, dx +\int_{\Omega}((\z_n(t)-\z_m(t)),D(u_{n}(t)-u_m(t)) \\
- \int_{\partial\Omega}(u_{n}(t)-u_m(t))[\z_n(t)-\z_m(t ),\nu]\, d\mathcal H^{N-1} =\int_{\Omega}(f_n(t)-f_m(t))(u_{n}(t)-u_m(t))\, dx\,.
\end{multline*}
and Proposition \ref{Prop.2.9} implies
\begin{multline*}
\frac{1}{2}\bigg(\int_{\Omega}(u_{n}(t)-u_m(t))^{2}\, dx\bigg)^\prime +\int_{\Omega}((\z_n(t)-\z_m(t)),D(u_{n}(t)-u_m(t)) \\
- \int_{\partial\Omega}(u_{n}(t)-u_m(t))[\z_n(t)-\z_m(t ),\nu]\, d\mathcal H^{N-1} =\int_{\Omega}(f_n(t)-f_m(t))(u_{n}(t)-u_m(t))\, dx\,.
\end{multline*}
Integrating between $0$ and $t\in (0,T)$, dropping two nonnegative terms and having in mind that the initial data are the same, we obtain
\begin{multline*}
  \frac{1}{2}\int_{\Omega}(u_{n}(t)-u_m(t))^{2}\, dx\le
\int_0^t\int_\Omega|f_n(t)-f_m(t)||u_{n}(t)-u_m(t)|\, dx\, dt \\
  \le \int_0^T\int_\Omega|f_n(t)-f_m(t)||u_{n}(t)-u_m(t)|\, dx\, dt\,.
\end{multline*}
Now the right hand side tends to $0$ since $ f_{n}\to f $  in $ L^1(0,T;L^2(\Omega)) $ and $\{u_n\}$ is bounded in $ L^{\infty}(0,T;L^2(\Omega))$.  We conclude that $\{u_n\}_{n=1}^{\infty}$ is
a Cauchy sequence in $ L^{\infty}(0,T;L^2(\Omega))$ and so there exists $u\in C([0,T];L^2(\Omega))$ such that
\begin{equation}\label{conv_C}
  u_n\longrightarrow u \quad\text{ in } C([0,T]; L^2(\Omega))\,.
\end{equation}
  As a consequence, the function
$u(t)$ is well--defined for all $t\in[0,T]$.

\subsection{Step 3: Convergence $u_n^2 \to u^2$ in $L^\infty(0,T;  L^1(\Omega))$}

This technical consequence of the previous Step 2 will be used in Step 14.
Recalling \eqref{mainEstimate} and applying H\"older's inequality,
\begin{multline*}
\int_\Omega |u_n(t)^2-u(t)^2| \,dx = \int_\Omega |u_n(t)-u(t)||u_n(t)+u(t)|\,dx \\
\le \left(\int_\Omega |u_n(t)-u(t)|^2 \,dx\right)^\frac{1}{2}\left(\int_\Omega |u_n(t)+u(t)|^2\,dx\right)^{\frac12}\le 2C \left(\int_\Omega |u_n(t)-u(t)|^2 \,dx\right)^\frac{1}{2}.
\end{multline*}
Step 3 is now straightforward.

\subsection{Step 4: $u \in L_w^1(0,T;BV(\Omega))$}

Going back to \eqref{conv_C}, we deduce
\begin{equation*}
u_n(t)\longrightarrow u(t) \quad\text{ in } L^1(\Omega) \text{ for every } t\in(0,T)\,.
\end{equation*}
It follows from this convergence and the lower semicontinuity of the total variation that
\begin{equation*}
\int_\Omega|Du(t)|\le \liminf_{n\to\infty}\int_\Omega|Du_n(t)|
\end{equation*}
holds for almost all $t\in(0,T)$.
Applying Fatou's lemma, we get
\begin{equation*}
\int_0^T\int_\Omega |Du(t)|\,dt   \le \int_0^T \liminf_{n\to\infty} \int_\Omega |Du_n(t)|\, dt   \le \liminf_{n\to\infty} \int_0^T\int_\Omega |Du_n(t)| \, dt \le C_2\,.
\end{equation*}
Thus, $u(t)\in BV(\Omega)$ for almost all $t\in(0,T)$ and so \cite[Lemma 5.19]{ACM} implies the function $t\mapsto \int_\Omega |Du(t)|$ is measurable and $u\in L_w^1(0,T;BV(\Omega))$.

\subsection{Step 5: Existence of the vector field $\z\in L^\infty((0,T)\times\Omega;\R^N)$}

Our next objective is to see that equation holds in the sense of distributions. To this end, we need to get the vector field $\z\in L^\infty((0,T)\times\Omega;\R^N)$ which plays the role of $Du/|Du|$, and the element $\xi$ that plays the role of the time derivative of $u$. In addition, we establish the sense in which $\z_n$ converges to $\z$ and $u_n'$ converges to $\xi$. The easy work corresponding to $\z$ will be done in this Step, while the corresponding to $\xi$ in Steps 6-7. Finally, in Step 8, we check that the equation holds in the sense of distributions.

\bigskip

\noindent For every $n\in\N$, it holds $\|\z_n\|_{\infty}\le 1$ then, up to a subsequence, $\z_n \estrella\z$ in $L^\infty((0,T)\times\Omega;\R^N)$ and $\|\z\|_\infty\le 1$.

\subsection{Step 6:  Convergence $\Div\z_\alpha \estrella \Div\z$ in $L^1(0,T;BV(\Omega)\cap L^2(\Omega))^*$ for some subnet $\{\z_\alpha\}_{\alpha\in I}$}

Let $v\in L^1(0,T;BV(\Omega)\cap L^2(\Omega))$. Since for almost every $t\in(0,T)$, $\Div\z_n(t)\in L^2(\Omega)$, it follows from Anzellotti's theory that
\begin{equation}\label{ec3}
-\int_\Omega v(t) \Div\z_n(t)\,dx=
\int_\Omega (\z_n(t),Dv(t))  - \int_{\partial\Omega} v(t) [\z_n(t),\nu]\,d\h
\end{equation}
for almost every $t\in(0,T)$. Then, integrating between $0$ and $T$, it becomes
\begin{align*}
\left| \int_0^T \int_\Omega v(t) \Div\z_n(t)\,dx,\,dt \right| & \le \int_0^T \left\{\int_\Omega |(\z_n(t),Dv(t))|
+ \int_{\partial\Omega} |v(t)||[\z_n(t),\nu]|\,d\h \right\}\,dt \\
& \le  \|\z_n\|_{\infty}\int_0^T\left\{\int_\Omega |Dv(t)|
+ \int_{\partial\Omega} |v(t)|\,d\h \right\}\,dt \\
& \le \int_0^T \|v(t)\|_{BV(\Omega)}\, dt\le \int_0^T \|v(t)\|_{BV(\Omega)\cap L^2(\Omega)}\, dt\,,
\end{align*}
and so the sequence $\{\Div\z_n\}_{n=1}^\infty$ is bounded in $ L^1(0,T;BV(\Omega)\cap L^2(\Omega))^*$.
Therefore, there exists $\rho\in L^1(0,T;BV(\Omega)\cap L^2(\Omega))^*$ and a subnet such that
\begin{equation*}
\Div \z_\alpha\estrella\rho \quad\text{ in } L^1(0,T;BV(\Omega)\cap L^2(\Omega))^*\,.
\end{equation*}
Let $\omega\in C_0^\infty(\Omega)$ and $\eta\in C_0^\infty(0,T)$. Due to \eqref{ec3}, we already know that
\begin{equation*}
\int_0^T\int_\Omega\eta(t)\,\omega\,\Div\z_\alpha(t) \,dt = -\int_0^T \eta(t)\int_\Omega\z_\alpha(t)\cdot\nabla \omega \,dx\,dt\,,
\end{equation*}
and taking limits in $\alpha\in I$ the equality becomes
\begin{equation*}
\int_0^T\eta(t)\langle\,\rho(t),\omega\,\rangle_\Omega \,dt = -\int_0^T\eta(t)\int_\Omega\z(t)\cdot\nabla \omega \,dx\,dt
\end{equation*}
which implies $\langle\,\rho(t),\omega\,\rangle_\Omega=\langle\,\Div\z(t),\omega\, \rangle_\Omega$ for almost every $t\in(0,T)$ and every $\omega\in C_0^\infty(\Omega)$.\\
 Observe that $\Div\z$ can be extended uniquely to an element of $L^1 (0,T; W_0^{1,1}(\Omega)\cap L^2(\Omega))^*$ and we have proven that $\rho$ is one of the further extensions of $\Div\z$ to the bigger space $ L^1 (0,T; BV(\Omega)\cap L^2(\Omega))^*$. Since this extension not need to be unique, from now on, we will identify $\Div\z$  with this specific extension.

\bigskip

\noindent Now, we define the element which performs the role of the time derivative of $u$:
\begin{equation}\label{Def_xi}
\xi = \Div\z+f \in L^1(0,T;BV(\Omega)\cap L^2(\Omega))^*+L^1(0,T;L^2(\Omega))\,.
\end{equation}

\subsection{Step 7: Convergence $u_\alpha^\prime \estrella \xi$ in $L^1(0,T;BV(\Omega)\cap L^2(\Omega))^*+L^1(0,T;L^2(\Omega))$}

 We start taking a test function $v\in L^1(0,T;BV(\Omega)\cap L^2(\Omega))\cap L^\infty(0,T;L^2(\Omega))$ in the equation $u_\alpha^\prime = \Div\z_\alpha + f_\alpha$ to obtain
\begin{equation*}
\int_0^T\int_\Omega u_\alpha^\prime(t)v(t)\,dx\,dt = \int_0^T\int_\Omega v(t)\Div\z_\alpha(t)\,dx\,dt + \int_0^T\int_\Omega f_\alpha(t)v(t)\,dx\,dt\,.
\end{equation*}
Taking limits in $\alpha\in I$ and considering convergences $\Div\z_\alpha \estrella \Div\z$ in the dual of $L^1(0,T;BV(\Omega)\cap L^2(\Omega))$ and $f_\alpha\to f$ in $L^1(0,T;L^2(\Omega))$ we get the desired result:
\begin{multline*}
\lim_{\alpha\in I} \int_0^T\int_\Omega u_\alpha^\prime(t)v(t)\,dx\,dt  = \int_0^T\langle\,\Div\z(t),v(t)\,\rangle_\Omega\,dt + \int_0^T\int_\Omega f(t)v(t)\,dx\,dt\\
=\int_0^T\ldoble \xi(t), v(t)\rdoble_\Omega\, dt\,.
\end{multline*}

\subsection{Step 8: For almost all $t$ the equation holds in the distributional sense}

 Let now $\omega\in C_0^\infty(\Omega)$ and $\eta\in C_0^\infty(0,T)$. We take the test function $\eta(t)\omega$ in $u_\alpha^\prime = \Div\z_\alpha(t) + f_\alpha$ to get
\begin{align*}
-\int_0^T \eta'(t)\int_\Omega u_\alpha(t) \omega \,dx \,dt & = \int_0^T\eta(t)\int_\Omega u_\alpha^\prime(t)\omega \,dx\,dt\\
&= \int_0^T\eta(t) \int_\Omega\z_\alpha(t)\cdot\nabla \omega \,dx\,dt + \int_0^T\eta(t)\int_\Omega f_\alpha(t)\,\omega \,dx\,dt\,.
\end{align*}
Since $u_\alpha \to u$ in $L^1((0,T)\times\Omega)$, $u_\alpha^\prime \estrella \xi$ in $L^1(0,T;BV(\Omega)\cap L^2(\Omega))^*+L^1(0,T;L^2(\Omega))$, $\z_\alpha \estrella\z$ in $L^\infty((0,T)\times\Omega;\R^N)$ and $f_\alpha\to f$ in $L^1(0,T;L^2(\Omega))$, taking limits in $\alpha\in I$ we arrive at
\begin{align*}
-\int_0^T \eta'(t)\int_\Omega u(t) \omega \,dx \,dt & = \int_0^T\eta(t)\ldoble\,\xi(t),\omega\,\rdoble_\Omega\,dt\\
&= \int_0^T\eta(t) \int_\Omega\z(t)\cdot\nabla \omega \,dx\,dt + \int_0^T\eta(t)\int_\Omega f(t)\,\omega \,dx\,dt\,
\end{align*}
for all $\eta\in C_0^\infty(0,T)$, which implies
\begin{equation*}
\left(\int_\Omega u(t) \omega \,dx \right)^\prime=\ldoble\,\xi(t),\omega\,\rdoble_\Omega= \int_\Omega\z(t)\cdot\nabla\omega \,dx + \int_\Omega f(t)\,\omega \,dx
\end{equation*}
for almost every $t\in(0,T)$.

\subsection{Step 9: $(\z(t),Dv)$ is a Radon measure in $\Omega$ for a.e. $t\in(0,T)$ and for all $v\in BV(\Omega)\cap L^2(\Omega)$}

The actual aim of this Step (and the following three ones) is to check that a Green's formula is available for  $\z(t)$. We point out that these vectors fields  does not satisfy the assumptions of \cite{An} since we cannot assure that its divergence is a Radon measure.

\noindent
Fix $v\in BV(\Omega)\cap L^2(\Omega) $, and consider $\omega\in C_0^\infty(\Omega)$ and $\eta\in C_0^\infty(0,T)$ with $\eta\ge0$. Recall that  $\Div\z_\alpha(t)\in L^2(\Omega)$ for almost all $t\in (0,T)$ and for all $\alpha$; consequently, due to Anzellotti's theory, $(\z_\alpha (t), Dv)$ is a Radon measure which satisfies
\begin{equation*}
\left|\int_\Omega \omega\,(\z_\alpha(t),Dv)\right| \le \|\omega\|_{L^\infty(\Omega)} \int_\Omega |Dv|\,
\end{equation*}
for almost every $t\in(0,T)$ (recall that $\|\z_\alpha\|_\infty\le 1$). Moreover, the following Green's formula holds
\begin{equation*}
-\langle\,\Div\z_\alpha(t),v\,\omega\,\rangle_\Omega - \int_\Omega v\,\z_\alpha(t)\cdot\nabla \omega\,dx = \int_\Omega \omega\,(\z_\alpha(t),Dv)
\end{equation*}
which implies
\begin{multline*}
\left| \int_0^T\eta(t)\Big[ \langle\,\Div\z_\alpha(t),v\,\omega\,\rangle_\Omega + \int_\Omega v\,\z_\alpha(t)\cdot\nabla \omega\,dx \Big] \,dt \right| \le \int_0^T\eta(t)\left|  \int_\Omega \omega\,(\z_\alpha(t),Dv)\right|\,dt \\
 \le \|\omega\|_{L^\infty(\Omega)} \int_0^T \eta(t)\int_\Omega |Dv|\,dt < +\infty\,.
\end{multline*}

Now, we take limits in $\alpha\in I$ to get
\begin{equation*}
\left| \int_0^T\eta(t)\Big[ \langle\,\Div\z(t),v\,\omega\,\rangle_\Omega + \int_\Omega v\,\z(t)\cdot\nabla \omega\,dx \Big] \,dt \right| \le \|\omega\|_{L^\infty(\Omega)}\int_0^T \eta(t) \int_\Omega |Dv|\,dt\,.
\end{equation*}
We deduce that for almost every $t\in(0,T)$ it holds:
\begin{equation*}
\left| \langle\,\Div\z(t),v\,\omega\,\rangle_\Omega + \int_\Omega v\,\z(t)\cdot\nabla \omega\,dx \right| \le \|\omega\|_{L^\infty(\Omega)}\int_\Omega |Dv|\,,
\end{equation*}
 from where it follows
\begin{equation*}
\Big| \int_\Omega \omega\,(\z (t),Dv)\Big| \le \|\omega\|_{L^\infty(\Omega)}\int_\Omega |Dv|
\end{equation*}
and so $(\z(t),Dv)$ is a Radon measure in $\Omega$.

\subsection{Step 10: Definition of the trace on the boundary of the normal component $[\z(t), \nu]$}

For every $v\in W^{1,1}(\Omega)\cap L^\infty(\Omega)$, we define
\begin{equation}\label{def_frontera}
\langle\,\z(t), v\,\rangle_{\partial\Omega}=\langle\, \Div\z(t), v\,\rangle_\Omega+\int_\Omega\z(t)\cdot\nabla v\, dx\,.
\end{equation}
We point out that this value also depends on the extension $\rho$, which we have identified with $\Div\z$.

\noindent
Now, let $\eta\in C_0^\infty (0, T)$ be nonnegative and  consider
\begin{equation*}
\int_0^T\eta(t)\langle\,\z(t), v\,\rangle_{\partial\Omega}\,dt=\int_0^T\eta(t)\langle\, \Div\z(t), v\,\rangle_\Omega\, dt+\int_0^T\eta(t)\int_\Omega\z(t)\cdot\nabla v\, dx\, dt\,.
\end{equation*}
Notice that if $v_1, v_2\in W^{1,1}(\Omega)\cap L^\infty(\Omega)$ satisfy $v_1=v_2$ on $\partial\Omega$, then
\begin{equation*}
\int_0^T\eta(t)\langle\, \Div\z(t), v_1-v_2\,\rangle_\Omega\, dt+\int_0^T\eta(t)\int_\Omega\z(t)\cdot\nabla (v_1-v_2)\, dx\, dt =0
\end{equation*}
since we take $\Div\z(t)$ in the distributional sense and $v_1-v_2\in W_0^{1,1}(\Omega)$.
Therefore,
\begin{equation*}
\int_0^T\eta(t)\langle\,\z(t), v_1\,\rangle_{\partial\Omega}\,dt=\int_0^T\eta(t)\langle\,\z(t), v_2\,\rangle_{\partial\Omega}\,dt
\end{equation*}
and so $\displaystyle\int_0^T\eta(t)\langle\,\z(t), v\,\rangle_{\partial\Omega}\,dt$ only depends on $v$ through its trace.

\noindent
On the other hand, given $\alpha\in I$, since $\Div\z_\alpha(t)\in L^2(\Omega)$ a.e., Anzellotti's theory applies and so
\begin{equation*}
|\langle\,\z_\alpha(t), v\,\rangle_{\partial\Omega}|\le \|\z_\alpha(t)\|_\infty \int_{\partial\Omega}|v|\, d\mathcal H^{N-1}\le \int_{\partial\Omega}|v|\, d\mathcal H^{N-1}
\end{equation*}
wherewith
\begin{multline*}
\left|\int_0^T\eta(t)\langle\, \Div\z_\alpha(t), v\,\rangle_\Omega\, dt+\int_0^T\eta(t)\int_\Omega\z_\alpha(t)\cdot\nabla v\, dx\, dt\right| \\
\le \int_0^T \eta(t)\int_{\partial\Omega}|v|\, d\mathcal H^{N-1}\, dt<+\infty\,.
\end{multline*}
Taking the limit for $\alpha\in I$, it yields
\begin{equation*}
\left|\int_0^T\eta(t)\langle\, \Div\z(t), v\,\rangle_\Omega\, dt+\int_0^T\eta(t)\int_\Omega\z(t)\cdot\nabla v\, dx\, dt\right|\le
\int_0^T \eta(t)\int_{\partial\Omega}|v|\, d\mathcal H^{N-1}\, dt\,,
\end{equation*}
and consequently we deduce that
\begin{equation*}
\left|\int_0^T\eta(t)\langle\,\z(t), v\,\rangle_{\partial\Omega}\, dt\right|\le \int_0^T\eta(t)\int_{\partial\Omega}|v|\, d\mathcal H^{N-1}dt
\end{equation*}
for every nonnegative test function $\eta\in C_0^\infty(\Omega)$.

\noindent
In the same spirit of \cite{An}, for each $t\in(0,T)$, we define $F_t\colon L^\infty(\partial\Omega)\to\R$ by
\begin{equation*}
F_t(w)=\langle\,\z(t), v\,\rangle_{\partial\Omega}
\end{equation*}
where $w\in L^\infty(\partial\Omega)$ and $v\in W^{1,1}(\Omega)\cap L^\infty(\Omega)$ satisfies $v\big|_{\partial\Omega}=w$. We have seen that
\begin{equation*}
\int_0^T\eta(t) F_t(w)\, dt
\end{equation*}
is well--defined and
\begin{equation*}
\left|\int_0^T\eta(t) F_t(w)\, dt\right|\le \int_0^T\eta(t)\int_{\partial\Omega}|w|\, d\mathcal H^{N-1}dt
\end{equation*}
for all $\eta\in C_0^\infty(0,T)$.
Hence, for every $w\in L^\infty(\partial\Omega)$,
\begin{equation*}
|F_t(w)|\le \int_{\partial\Omega}|w|\, d\mathcal H^{N-1}
\end{equation*}
holds for almost all $t\in (0,T)$. Note that the null set depends on $w$.

To go on we have to use a separability argument. Let $V$ denote a countable set which is dense in $W^{1,1}(\Omega)$. Truncating functions of $V$, if necessary, we may assume that $v\in W^{1,1}(\Omega)\cap L^\infty(\Omega)$ for all $v\in V$. We now get for almost all $t\in (0,T)$:
\begin{equation*}
|F_t(w)|\le \int_{\partial\Omega}|w|\, d\mathcal H^{N-1}
\end{equation*}
for every $w\in L^\infty(\partial\Omega)$ satisfying $w=v\big|_{\partial\Omega}$ with $v\in V$. Next fix one of these points $t$, choose $w_0\in L^\infty(\partial\Omega)$ and let $v_0\in W^{1,1}(\Omega)\cap L^\infty(\Omega)$ such that its trace is $w_0$. Consider a sequence $(v_n)_n$ in $V$ satisfying $v_n\to v_0$ in $W^{1,1}(\Omega)$. It leads to
\[\begin{array}{ll}
& \langle \Div\z(t), v_n\rangle_\Omega\to \langle \Div\z(t), v_0\rangle_\Omega\\
& \nabla v_n\to \nabla v_0 \hbox{ in } L^1(\Omega; \R^n)
\end{array}\]
 wherewith $F_t(w_n)\to F_t(w_0)$. On the other hand, $v_n\to v_0$ in $W^{1,1}(\Omega)$ also implies $w_n\to w_0$ in $L^{1}(\partial\Omega)$. Thus, it follows from $|F_t(w_n)|\le \int_{\partial\Omega}|w_n|\, d\mathcal H^{N-1}$
for all $n\in \N$ that
$|F_t(w)|\le \int_{\partial\Omega}|w|\, d\mathcal H^{N-1}$.
Therefore,
\[|F_t(w)|\le \int_{\partial\Omega}|w|\, d\mathcal H^{N-1}\qquad \forall w\in L^\infty(\partial\Omega)\]
holds for almost all $t\in (0,T)$.
Taking one of these $t\in (0,T)$,
the functional $F_t$ may be extended to a functional in $L^1(\partial\Omega)^*$, so that is represented by a $L^\infty$--function, denoted by $[\z(t),\nu]$. In other words, $[\z(t), \nu]\in L^\infty(\partial\Omega)$ in such a way that $\|[\z(t), \nu]\|_\infty\le 1$ and
\begin{equation*}
F_t(w)=\int_{\partial\Omega}w [\z(t),\nu]\, d\mathcal H^{N-1}\qquad \forall w\in L^\infty(\partial\Omega)
\end{equation*}
for almost all $t\in (0,T)$. Moreover, we have deduced the following Green's formula holds  for almost all $t\in (0,T)$:
\begin{equation}\label{ec:8.1}
\langle\,\Div\z(t),v\,\rangle_\Omega+\int_\Omega\z(t)\cdot\nabla v \,dx =  \int_{\partial\Omega}v[\z(t),\nu]\,d\h
\end{equation}
for every $v\in W^{1,1}(\Omega)\cap L^\infty(\Omega)$.

\subsection{Step 11: Convergence of the traces on the boundary of the normal components}

Let $\eta \in C_0^\infty (0,T)$ and $w\in L^\infty(\partial\Omega)$. We will prove that
\begin{equation}\label{Conv_frontera}
\lim_{\alpha\in I} \int_0^T\eta(t)\int_{\partial \Omega} w[\z_\alpha(t),\nu]\,d\h \,dt =\int_0^T\eta(t) \int_{\partial \Omega} w\,[\z(t),\nu]\,d\h \,dt\,.
\end{equation}
If $v\in W^{1,1}(\Omega)\cap L^\infty(\Omega)$ is such that $v\big|_{\partial\Omega}=w$, then for every $\alpha\in I$ the Green's formula holds
\begin{multline*}
\int_0^T \eta(t)\int_\Omega\z_\alpha(t)\cdot\nabla v \,dx\,dt + \int_0^T\eta(t)
\int_\Omega v\,\Div\z_\alpha(t)\,dx\,dt\\
 = \int_0^T\eta(t)\int_{\partial\Omega}w[\z_\alpha(t),\nu]\,d\h\,dt\,.
\end{multline*}
We can take limits in $\alpha\in I$ on the left hand side to get
\begin{align*}
\lim_{\alpha\in I} \int_0^T& \eta(t)\int_{\partial\Omega}w[\z_\alpha(t),\nu]\,d\h\,dt \\
&=\lim_{\alpha\in I} \left[\int_0^T \eta(t)\int_\Omega\z_\alpha(t)\cdot\nabla v \,dx\,dt + \int_0^T\eta(t)\int_\Omega v\,\Div\z_\alpha(t)\,dx\,dt \right] \\
&= \int_0^T \eta(t)\int_\Omega\z(t)\cdot\nabla v \,dx\,dt + \int_0^T\eta(t)\langle\,\Div\z(t),v\,\rangle_\Omega\,dt  \\
&=  \int_0^T\eta(t)\int_{\partial\Omega}v[\z(t),\nu]\,d\h\,dt=  \int_0^T\eta(t)\int_{\partial\Omega}w[\z(t),\nu]\,d\h\,dt\,,
\end{align*}
where we have used Green's formula \eqref{ec:8.1}.

\subsection{Step 12: Green's formula }

Let $v\in BV(\Omega)\cap L^2(\Omega)$. We are going to show that
\begin{equation*}
\langle\, \Div\z(t),v\,\rangle_\Omega + \int_\Omega (\z(t),Dv) = \int_{\partial\Omega} v[\z(t),\nu]\,d\h
\end{equation*}
holds for almost every $t\in(0,T)$.
\medskip

\noindent Consider $\eta\in C_0^\infty(0,T)$. Owing to \eqref{ec3}, for $k>0$ and for every $\alpha\in I$, it holds
\begin{multline}\label{nivelA}
\int_0^T\eta(t)
\int_\Omega T_k(v)\,\Div\z_\alpha(t) \, dx\,dt + \int_0^T\eta(t)\int_\Omega (\z_\alpha(t),DT_k(v))\,dt \\
= \int_0^T\eta(t)\int_{\partial\Omega} T_k(v)[\z_\alpha(t),\nu]\,d\h\,dt\,.
\end{multline}
We remark the needed to use $T_k(v)$ instead of $v$ in order to handle the integral on the boundary.

\noindent
Our aim is to take the limit in $\alpha\in I$ in this identity. The limit in the first term is consequence of $\Div\z_\alpha\estrella \Div\z$ in $L^1(0,T;BV(\Omega)\cap L^2(\Omega))^*$. To deal with the second term we may argue as in \cite[Proposition 2.1]{An} since
\begin{equation*}
\int_0^T\eta(t)\int_U |(\z_\alpha(t), DT_k(v))|\, dt\le \int_0^T\eta(t)\int_U|DT_k(v)|\, dt
\end{equation*}
for all open $U\subset\Omega$ and
\begin{equation*}
\lim_{\alpha\in I} \int_0^T\eta(t)\int_\Omega \omega(\z_\alpha(t),DT_k(v))\, dt= \int_0^T\eta(t)\int_\Omega \omega (\z(t),DT_k(v))\, dt
\end{equation*}
for all $\omega\in C_0^\infty(\Omega)$. Finally, on the right hand side, we may apply  \eqref{Conv_frontera} and so  we are able to take limits in $\alpha\in I$, wherewith \eqref{nivelA} becomes
\begin{multline}\label{nivelK}
\int_0^T\eta(t)\langle\,\Div\z(t),T_k(v)\,\rangle_\Omega \,dt + \int_0^T\eta(t)\int_\Omega (\z(t),DT_k(v))\,dt \\
= \int_0^T\eta(t)\int_{\partial\Omega} T_k(v)[\z(t),\nu]\,d\h\,dt\,.
\end{multline}
Having in mind that $T_k(v)\to v$ in $BV(\Omega)\cap L^2(\Omega)$, we may let
 $k$ go to $\infty$ in \eqref{nivelK}. Thus,
\begin{equation*}
\int_0^T\eta(t)\langle\,\Div\z(t),v\,\rangle_\Omega \,dt + \int_0^T\eta(t)\int_\Omega (\z(t),Dv)\,dt = \int_0^T\eta(t)\int_{\partial\Omega} v[\z(t),\nu]\,d\h\,dt\,
\end{equation*}
for all $\eta\in C_0^\infty(0,T)$, which implies
\begin{equation*}
\langle\,\Div\z(t),v\,\rangle_\Omega  + \int_\Omega (\z(t),Dv) =\int_{\partial\Omega} v[\z(t),\nu]\,d\h
\end{equation*}
for almost all $t\in(0,T)$.

\subsection{Step 13: $\xi$ is the time derivative of $u$ in the sense of Definition \ref{T.Der}}

Let $\Psi\in L^1(0,T;BV(\Omega)\cap L^2(\Omega))\cap L^\infty(0,T;L^2(\Omega))$ with compact support in $(0,T)$ and let $\Theta\in L^1_w(0,T;BV(\Omega))\cap L^1(0,T;L^2(\Omega))$ be the weak derivative of $\Psi$. Since
\begin{equation*}
\int_0^T\int_\Omega u_\alpha^\prime(t) \Psi(t)\,dx\,dt = -\int_0^T\int_\Omega u_\alpha(t) \Theta(t)\,dx\,dt\,,
\end{equation*}
 $u_\alpha^\prime \estrella \xi$ in $L^1(0,T;BV(\Omega)\cap L^2(\Omega))^*+L^1(0,T;L^2(\Omega))$ and $u_\alpha \to u$ in $L^\infty(0,T; L^2(\Omega))$, we can take limits in $\alpha\in I$ to obtain
\begin{equation*}
\int_0^T \ldoble \xi(t), \Psi(t) \rdoble_\Omega\,dt = -\int_0^T\int_\Omega u(t) \,\Theta(t)\,dx\,dt\,.
\end{equation*}

\subsection{Step 14: Conditions \eqref{cond2} and \eqref{cond3}}

To prove that $u$ is a solution to problem \eqref{prob}, it remains to check that satisfies conditions \eqref{cond2} and \eqref{cond3}  of Definition \ref{def0}.

Let $\eta\in C_0^\infty(0,T)$ be nonnegative and define
\begin{equation*}
\Psi_\varepsilon(t)=\frac1{\varepsilon}\int_{t-\varepsilon}^t\eta(s) u(s)\, ds\,,
\end{equation*}
that is a Pettis integral.
Reasoning as in the beginning of Proposition \ref{para_unic}, we see that
\begin{multline}\label{ec:15.0}
-\frac1{\varepsilon}\int_0^T\big(\eta(t-\varepsilon)-\eta(t)\big)\int_\Omega  \frac12 u(t)^2\, dx\, dt\le -\int_0^T\ldoble \,\xi(t), \Psi_\varepsilon(t)\,\rdoble_\Omega\, dt\\
=- \int_0^T  \frac1{\varepsilon}\int_{t-\varepsilon}^t\eta(s)\Big[ \langle \Div \z(t), u(s)\rangle_\Omega+ \int_\Omega f(t)u(s)\, dx\Big]\, ds\, dt\\
=\int_0^T  \frac1{\varepsilon}\int_{t-\varepsilon}^t\eta(s) \Big[\int_\Omega (\z(t), Du(s))-\int_{\partial\Omega}u(s)[\z(t),\nu]\, d\mathcal H^{N-1}-\int_\Omega u(s) f(t)\, dx\Big]\, ds\, dt\\
\le  \int_0^T  \frac1{\varepsilon}\int_{t-\varepsilon}^t\eta(s) \Big[\int_\Omega |Du(s)|-\int_{\partial\Omega}u(s)[\z(t),\nu]\, d\mathcal H^{N-1}-\int_\Omega u(s) f(t)\, dx\Big]\, ds\, dt\,.
\end{multline}
Having in mind Lemma \ref{lema0} and Corollary \ref{corollary1}, we let $\varepsilon$ go to 0 to obtain
\begin{multline}\label{ec:15.1}
\int_0^T \eta'(t)\int_\Omega \frac12 u(t)^2\,dx\, dt\\
\le \int_0^T  \eta(t) \Big[\int_\Omega |Du(t)|-\int_{\partial\Omega}u(t)[\z(t),\nu]\, d\mathcal H^{N-1}-\int_\Omega u(t) f(t)\, dx\Big]\, dt\\
\le \int_0^T  \eta(t) \Big[\int_\Omega |Du(t)|+\int_{\partial\Omega}|u(t)|\, d\mathcal H^{N-1}-\int_\Omega u(t) f(t)\, dx\Big]\, dt\,.
\end{multline}
On the other hand, taking $\eta(t)u_n(t)$ as test function in problem \eqref{Pn}, it yields
\begin{multline*}
-\int_0^T \eta'(t)\int_\Omega \frac12 u_n(t)^2\, dx\, dt+\int_0^T \eta(t) \Big[\int_\Omega |Du_n(t)|+\int_{\partial\Omega}|u_n(t)|\, d\mathcal H^{N-1}\Big]\,dt\\
=\int_0^T\eta(t)\int_\Omega u_n(t)f_n(t)\, dx\, dt\,.
\end{multline*}
Applying Step 3, the lower semicontinuity of the BV-norm and Fatou's lemma, we deduce that
\begin{multline*}
-\int_0^T \eta'(t)\int_\Omega \frac12 u(t)^2\, dx\, dt+\int_0^T \eta(t) \left[\int_\Omega |Du(t)|+\int_{\partial\Omega}|u(t)|\, d\mathcal H^{N-1}\right]dt\\
\le \int_0^T\eta(t)\int_\Omega u(t) f(t)\, dx\, dt\,.
\end{multline*}
Finally, having in mind \eqref{ec:15.1}, it implies
\begin{align}\label{ec:15.2}
\int_0^T \eta'(t)&\int_\Omega \frac12 u(t)^2\,dx\, dt\\
\notag &\le \int_0^T  \eta(t) \Big[\int_\Omega |Du(t)|-\int_{\partial\Omega}u(t)[\z(t),\nu]\, d\mathcal H^{N-1}-\int_\Omega u(t) f(t)\, dx\Big]\, dt\\
\notag &\le \int_0^T  \eta(t) \Big[\int_\Omega |Du(t)|+\int_{\partial\Omega}|u(t)|\, d\mathcal H^{N-1}-\int_\Omega u(t) f(t)\, dx\Big]\, dt\\
\notag &\le\int_0^T \eta'(t)\int_\Omega \frac12 u(t)^2\,dx\, dt\,.
\end{align}

It follows from \eqref{ec:15.2} that
\begin{equation*}
-\int_0^T  \eta(t) \int_{\partial\Omega}u(t)[\z(t),\nu]\, d\mathcal H^{N-1}\, dt=\int_0^T  \eta(t) \int_{\partial\Omega}|u(t)|\, d\mathcal H^{N-1}\, dt\,.
\end{equation*}
Since this identity holds for every nonnegative $\eta\in C_0^\infty(\Omega)$, we get
\begin{equation*}
\int_{\partial\Omega}\big(|u(t)|+u(t)[\z(t),\nu]\big)\, d\mathcal H^{N-1}=0
\end{equation*}
for almost all $t\in(0,T)$, which implies the boundary condition \eqref{cond3}.

Another consequence of \eqref{ec:15.2} is the identity
\begin{multline*}
\int_0^T \eta'(t)\int_\Omega \frac12 u(t)^2\,dx\, dt\\
=\int_0^T  \eta(t) \Big[\int_\Omega |Du(t)|+\int_{\partial\Omega}|u(t)|\, d\mathcal H^{N-1}-\int_\Omega u(t) f(t)\, dx\Big]\, dt\,.
\end{multline*}
Notice that the arbitrariness of $\eta$ leads to condition \eqref{cond2} and so
Theorem \ref{exis} is now completely proven.
\end{pf}

\section{Uniqueness of solution}

In this section we show the uniqueness of the solution to problem \eqref{prob}.
\begin{Theorem}
For every $f\in L^1(0,T;L^2(\Omega))$ and every $u_0\in L^2(\Omega)$, there exists at most a solution to problem \eqref{prob}.
\end{Theorem}

\begin{pf}
Assume that $u_1$ and $u_2$ are two solutions to problem \eqref{prob}. Then there exist $\xi_1$ and $\xi_2$ which are the time derivatives of $u_1$ and $u_2$, respectively, and there also exist the corresponding vector fields $\z_1$ and $\z_2$. We also point out that $u_1(0)=u_2(0)$.
To see that $u_1=u_2$,  we fix $\eta\in C_0^\infty(0,T)$ such that $\eta \ge0$.
The proof is split into several stages.

\medskip
\textbf{Step 1: }
First we choose $\varepsilon>0$ so small for the following calculations to be held and define
\begin{equation*}
\Psi^1_\varepsilon(t)=\frac1{\varepsilon}\int_{t-\varepsilon}^t\eta(s) u_1(s) \, ds\,,\qquad \Psi^2_\varepsilon(t)=\frac1{\varepsilon}\int_{t-\varepsilon}^t\eta(s) u_2(s) \, ds\,,
\end{equation*}
which are actually Pettis integrals.

\noindent
Since $\xi_1$ and $\xi_2$ are the time derivative of $u_1$ and $u_2$, respectively, and $u_1$ and $u_2$ are solutions, it follows from \eqref{unic:0} that
\begin{equation}\label{unic:1}
\lim_{\varepsilon\to0}\int_0^T\ldoble\, \xi_1(t), \Psi^1_\varepsilon (t) \,\rdoble_\Omega\, dt=
-\int_0^T\eta'(t)\int_\Omega  \frac12 u_1(t)^2\, dx\, dt
\end{equation}
and
\begin{equation}\label{unic:2}
\lim_{\varepsilon\to0}\int_0^T\ldoble\, \xi_2(t), \Psi^2_\varepsilon (t) \,\rdoble_\Omega\, dt=
-\int_0^T\eta'(t)\int_\Omega  \frac12 u_2(t)^2\, dx\, dt\,.
\end{equation}
On the other hand, we have
\begin{multline*}
\int_0^T\ldoble\, \xi_1(t)+\xi_2(t), \Psi^1_\varepsilon (t)+\Psi^2_\varepsilon(t) \,\rdoble_\Omega\, dt\\
=-\frac1{\varepsilon}
\int_0^T \int_\Omega \Big(\eta(t)\big(u_1(t)+u_2(t)\big)-\eta(t-\varepsilon) \big(u_1(t-\varepsilon)+u_2(t-\varepsilon)\big)\Big) \big(u_1(t)+u_2(t)\big)\, dx\, dt\\
=-\frac1{\varepsilon}
\int_0^T \int_\Omega \eta(t)\big(u_1(t)+u_2(t)\big)^2\, dx \, dt+\frac1{\varepsilon}
\int_0^T \int_\Omega\eta(t-\varepsilon) \big(u_1(t-\varepsilon)+u_2(t-\varepsilon)\big)\big( u_1(t)+u_2(t)\big)\, dx\, dt\\
=\frac1{\varepsilon} \int_0^T \eta(t) \int_\Omega \big(u_1(t)+u_2(t)\big)\Big(\big(u_1(t+\varepsilon)+u_2(t+\varepsilon)\big)-\big(u_1(t)+u_2(t)\big)\Big)\, dx\, dt\,.
\end{multline*}
Since $\big(u_1(t+\varepsilon)+u_2(t+\varepsilon)\big)\big(u_1(t)+u_2(t) \big)\le\frac{1}{2}\big(u_1(t+\varepsilon)+u_2(t+\varepsilon)\big)^2+\frac{1}{2}\big(u_1(t)+u_2(t)\big)^2$ holds, it yields
\begin{align*}
\int_0^T\ldoble\, \xi_1(t)&+\xi_2(t), \Psi^1_\varepsilon (t)+\Psi^2_\varepsilon(t) \,\rdoble_\Omega\, dt\\
&\le \frac1{\varepsilon} \int_0^T\int_\Omega \eta(t)  \Big(\frac12 \big(u_1(t+\varepsilon)+u_2(t+\varepsilon)\big)^2-\frac12 \big(u_1(t)+u_2(t)\big)^2\Big)\, dx\, dt\\
&=\frac1{\varepsilon} \int_0^T\int_\Omega  \frac12 \Big( \eta(t-\varepsilon)\big(u_1(t)+u_2(t)\big)^2-\eta(t) \big(u_1(t)+u_2(t)\big)^2\Big)\, dx\, dt\\
&=\frac1{\varepsilon} \int_0^T\big( \eta(t-\varepsilon)-\eta(t)\big)\int_\Omega  \frac12 \big(u_1(t)+u_2(t)\big)^2\, dx\, dt\\
&=\frac1{\varepsilon} \int_0^T\big( \eta(t-\varepsilon)-\eta(t)\big)\int_\Omega \Big[ \frac12 u_1(t)^2+\frac12 u_2(t)^2+u_1(t)u_2(t)\Big]\, dx\, dt\,.
\end{align*}
Letting $\varepsilon$ go to 0, we obtain
\begin{multline}\label{unic:3}
\limsup_{\varepsilon\to0}\int_0^T\ldoble\, \xi_1(t)+\xi_2(t), \Psi^1_\varepsilon (t)+\Psi^2_\varepsilon(t) \,\rdoble_\Omega\, dt\\
\le -\int_0^T\eta'(t)\int_\Omega \Big[ \frac12 u_1(t)^2+\frac12 u_2(t)^2+u_1(t)u_2(t)\Big]\, dx\, dt\,.
\end{multline}
Taking into account \eqref{unic:1} and \eqref{unic:2}, inequality \eqref{unic:3} becomes
\begin{equation}\label{unic:4}
\limsup_{\varepsilon\to0}\int_0^T\Big[\ldoble\, \xi_1(t), \Psi^2_\varepsilon(t) \,\rdoble_\Omega+\ldoble\, \xi_2(t), \Psi^1_\varepsilon(t) \,\rdoble_\Omega\Big]\, dt
\le -\int_0^T\!\eta'(t)\int_\Omega  u_1(t)u_2(t)\, dx\, dt\,.
\end{equation}

\noindent
Similarly,  we deduce

\begin{multline*}
\int_0^T\ldoble\, \xi_1(t)-\xi_2(t), \Psi^1_\varepsilon (t)-\Psi^2_\varepsilon(t) \,\rdoble_\Omega\, dt
\le \frac1{\varepsilon} \int_0^T\big( \eta(t-\varepsilon)-\eta(t)\big)\int_\Omega  \frac12 \big(u_1(t)-u_2(t)\big)^2\, dx\, dt\\
=\frac1{\varepsilon} \int_0^T\big( \eta(t-\varepsilon)-\eta(t)\big)\int_\Omega  \Big[\frac12 u_1(t)^2+\frac12 u_2(t)^2-u_1(t)u_2(t)\Big]\, dx\, dt
\end{multline*}
so that
\begin{multline*}
\limsup_{\varepsilon\to0}\int_0^T\ldoble\, \xi_1(t)-\xi_2(t), \Psi^1_\varepsilon (t)-\Psi^2_\varepsilon(t) \,\rdoble_\Omega\, dt\\
\le-\int_0^T\eta'(t)\int_\Omega \Big[ \frac12 u_1(t)^2+\frac12 u_2(t)^2-u_1(t)u_2(t)\Big]\, dx\, dt\,.
\end{multline*}
Then \eqref{unic:1} and \eqref{unic:2} imply
\begin{equation}\label{unic:5}
-\liminf_{\varepsilon\to0}\int_0^T\Big[\ldoble \,\xi_1(t), \Psi^2_\varepsilon(t) \rdoble_\Omega+\ldoble\, \xi_2(t), \Psi^1_\varepsilon(t) \,\rdoble_\Omega\Big]\, dt
\le \int_0^T\eta'(t)\int_\Omega  u_1(t)u_2(t)\, dx\, dt\,.
\end{equation}

\noindent
Gathering \eqref{unic:4} and \eqref{unic:5}, we conclude that there exists the limit and
\begin{equation*}
-\int_0^T\eta'(t)\int_\Omega  u_1(t)u_2(t)\, dx\, dt=
\lim_{\varepsilon\to0}\int_0^T\Big[\ldoble\, \xi_1(t), \Psi^2_\varepsilon(t) \,\rdoble_\Omega+\ldoble\, \xi_2(t), \Psi^1_\varepsilon(t) \,\rdoble_\Omega\Big]\, dt
\end{equation*}
which, together with \eqref{unic:1} and \eqref{unic:2}, turns out that
\begin{equation}\label{unic:6}
-\int_0^T\eta'(t)\int_\Omega  \frac12 \big( u_1(t)- u_2(t)\big)^2\, dx\, dt
=\lim_{\varepsilon\to0}\int_0^T\ldoble\, \xi_1(t)-\xi_2(t), \Psi^1_\varepsilon (t)-\Psi^2_\varepsilon(t) \,\rdoble_\Omega\, dt\,.
\end{equation}

\medskip
\noindent
\textbf{Step 2: }
Our next concern is to compute the limit, as $\varepsilon$ goes to 0, on the right hand side of \eqref{unic:6} in a different way. To this end, we write
\begin{multline}\label{unic:7}
\int_0^T\ldoble \,\xi_1(t)-\xi_2(t),\Psi_\varepsilon^1(t)-\Psi_\varepsilon^2(t)\,\rdoble_\Omega\, dt\\
=\int_0^T\ldoble \,\xi_1(t),\Psi_\varepsilon^1(t)\,\rdoble_\Omega\, dt+\int_0^T\ldoble \,\xi_2(t),\Psi_\varepsilon^2(t)\,\rdoble_\Omega\, dt\\
-\int_0^T\ldoble \,\xi_1(t),\Psi_\varepsilon^2(t)\,\rdoble_\Omega\, dt
-\int_0^T\ldoble \,\xi_2(t),\Psi_\varepsilon^1(t)\,\rdoble_\Omega\, dt=I_\varepsilon^a-I_\varepsilon^b
\end{multline}
where
\begin{equation*}
I_\varepsilon^a=\int_0^T\ldoble \,\xi_1(t),\Psi_\varepsilon^1(t)\,\rdoble_\Omega\, dt+\int_0^T\ldoble \,\xi_2(t),\Psi_\varepsilon^2(t)\,\rdoble_\Omega\, dt
\end{equation*}
and
\begin{equation*}
I_\varepsilon^b=\int_0^T\ldoble \,\xi_1(t),\Psi_\varepsilon^2(t)\,\rdoble_\Omega\, dt+\int_0^T\ldoble \,\xi_2(t),\Psi_\varepsilon^1(t)\,\rdoble_\Omega\, dt\,.
\end{equation*}
We now use the identities  $\xi_i(t)+\Div\z_i(t)=f(t)$ for $i=1,2$. With regard to $I_\varepsilon^a$, just apply again condition \eqref{unic:0} to get
\begin{multline}\label{unic:7.5}
\lim_{\varepsilon\to0}I_\varepsilon^a=-\int_0^T\eta'(t)\int_\Omega  \Big[\frac12 u_1(t)^2+\frac12 u_2(t)^2\Big]\, dx\, dt\\
=\int_0^T\eta(t)\Big[\ldoble \xi_1(t), u_1(t) \rdoble_\Omega+\ldoble \xi_2(t), u_2(t) \rdoble_\Omega\Big]\, dt\\
=\int_0^T\eta(t)\Big[\langle \Div\z_1(t), u_1(t) \rangle_\Omega+\langle \Div\z_2(t), u_2(t) \rangle_\Omega+\int_\Omega f(t) (u_1(t)+u_2(t))\, dx\Big]\, dt\\
=-\int_0^T\eta(t)\Big[\int_\Omega |Du_1(t)|+\int_\Omega |Du_2(t)|+\int_{\partial\Omega}|u_1(t)|\, d\mathcal H^{N-1}+\int_{\partial\Omega}|u_2(t)|\, d\mathcal H^{N-1}\Big]\, dt\\
+\int_0^T\eta(t)\int_\Omega f(t)(u_1(t)+u_2(t))\, dx\, dt\,.
\end{multline}
Notice that the existence of $\lim_{\varepsilon\to0}I_\varepsilon^b$ is now guaranteed by \eqref{unic:6} and \eqref{unic:7.5}, having in mind \eqref{unic:7}.

On the other hand, we have
\begin{equation*}
\int_0^T\ldoble \,\xi_1(t),\Psi_\varepsilon^2(t)\,\rdoble_\Omega\, dt=\int_0^T\frac1{\varepsilon}\int_{t-\varepsilon}^t\eta(s)\Big[\langle \Div\z_1(t), u_2(s)\rangle_\Omega+\int_\Omega f(t)u_2(s)\, dx\Big]\, ds\, dt\,.
\end{equation*}
Observing that
\begin{multline*}
-\int_0^T\frac1{\varepsilon}\int_{t-\varepsilon}^t\eta(s)\langle \Div\z_1(t), u_2(s)\rangle_\Omega\, ds\, dt\\
\le\int_0^T\frac1{\varepsilon}\int_{t-\varepsilon}^t\eta(s)\Big[\int_\Omega|(\z_1(t), Du_2(s)|+\int_{\partial\Omega}|u_2(s)[\z_1(t), \nu]|\, d\mathcal H^{N-1}\Big]\, ds\, dt\\
\le \int_0^T\frac1{\varepsilon}\int_{t-\varepsilon}^t\eta(s)\Big[\int_\Omega|Du_2(s)|+\int_{\partial\Omega}|u_2(s)|\, d\mathcal H^{N-1}\Big]\, ds\, dt\,,
\end{multline*}
we obtain that
\begin{multline*}
\int_0^T\ldoble \,\xi_1(t),\Psi_\varepsilon^2(t)\,\rdoble_\Omega\, dt\\
\ge -\int_0^T\frac1{\varepsilon}\int_{t-\varepsilon}^t\eta(s)\left[\int_\Omega|Du_2(s)|+\int_{\partial\Omega}|u_2(s)|\, d\mathcal H^{N-1}\right]\, ds\, dt\\
+\int_0^T\frac1{\varepsilon}\int_{t-\varepsilon}^t\eta(s)\int_\Omega f(t)u_2(s)\, dx\, ds\, dt
\end{multline*}
and, appealing to Lemma \ref{lema0}, the right hand side converges to
\begin{equation}\label{unic:8}
  -\int_0^T\eta(t)\left[\int_\Omega|Du_2(t)|+\int_{\partial\Omega}|u_2(t)|\, d\mathcal H^{N-1}\right]\, dt
+\int_0^T\eta(t)\int_\Omega f(t)u_2(t)\, dx\, dt
\end{equation}
as $\varepsilon$ goes to 0.
Analogously, we infer that
\begin{multline*}
\int_0^T\ldoble \,\xi_2(t),\Psi_\varepsilon^1(t)\,\rdoble_\Omega\, dt\\
\ge -\int_0^T\frac1{\varepsilon}\int_{t-\varepsilon}^t\eta(s)\left[\int_\Omega|Du_1(s)|+\int_{\partial\Omega}|u_1(s)|\, d\mathcal H^{N-1}\right]\, ds\, dt\\
+\int_0^T\frac1{\varepsilon}\int_{t-\varepsilon}^t\eta(s)\int_\Omega f(t)u_1(s)\, dx\, ds\, dt
\end{multline*}
whose right hand side converges, as $\varepsilon\to0$, to
\begin{equation}\label{unic:9}
-\int_0^T\eta(t)\left[\int_\Omega|Du_1(t)|+\int_{\partial\Omega}|u_1(t)|\, d\mathcal H^{N-1}\right]\, dt
+\int_0^T\eta(t)\int_\Omega f(t)u_1(t)\, dx\, dt\,.
\end{equation}
Finally, it follows from \eqref{unic:7.5}, \eqref{unic:8} and \eqref{unic:9} that
\begin{equation*}
\lim_{\varepsilon\to0}I_\varepsilon^b\ge  \lim_{\varepsilon\to0}I_\varepsilon^a\,.
\end{equation*}
Therefore, \eqref{unic:7} implies that
\begin{equation*}
\lim_{\varepsilon\to0}\int_0^T\ldoble \,\xi_1(t)-\xi_2(t),\Psi_\varepsilon^1(t)-\Psi_\varepsilon^2(t)\,\rdoble_\Omega\, dt= \lim_{\varepsilon\to0}I_\varepsilon^a- \lim_{\varepsilon\to0}I_\varepsilon^b\le 0\,.
\end{equation*}
Now, since \eqref{unic:6} holds,  it yields
\begin{multline*}
0\le \int_0^T\eta'(t)\int_\Omega \frac{1}{2}(u_1(t)-u_2(t))^2\,dx\,dt = -\int_0^T\eta(t)\left(\frac{1}{2}\int_\Omega (u_1(t)-u_2(t))^2\,dx\right)'\,dt \\= \lim_{\varepsilon\to0}I_\varepsilon^b- \lim_{\varepsilon\to0}I_\varepsilon^a
\end{multline*}
for all nonnegative $\eta\in C_0^\infty(0,T)$. Then,
\begin{multline*}
0\le -\left(\frac{1}{2}\int_\Omega (u_1(t)-u_2(t))^2\,dx\right)' \\
 \le 2\int_\Omega \left(|Du_1(t)| + |Du_2(t)|\right)+2\int_{\partial\Omega}\left( |u_1(t)|+|u_2(t)|\right)\,d\h \\
  = 2(\|u_1(t)\|_{BV(\Omega)} +\|u_2(t)\|_{BV(\Omega)})\in L^1(0,T)\,,
\end{multline*}
and we deduce that function $\displaystyle t\mapsto \int_\Omega  \big(u_1(t)-u_2(t)\big)^2\, dx $ is absolutely continuous in $(0,T)$.

\medskip
\noindent
\textbf{Conclusion:}
Since the function $t\mapsto \int_\Omega  \big(u_1(t)-u_2(t)\big)^2\, dx $ is absolutely continuous with nonpositive derivative, it follows that is nonincreasing in $(0,T)$.
Then
\begin{equation*}
\int_\Omega   \big( u_1(t)- u_2(t)\big)^2\, dx\, dt\le\int_\Omega  \big( u_1(0)- u_2(0)\big)^2\, dx\, dt=0\,\qquad \hbox{for all }t\in[0,T],
\end{equation*}
and we conclude that $u_1(t)=u_2(t)$ a.e. in $\Omega$ for every $t\in[0,T]$.
\end{pf}

As a consequence of this uniqueness result, each solution to problem \eqref{prob} can be obtained as a limit of solutions with $L^2$--data. It implies that every feature of these approximate solutions can be transferred to a general solution to  \eqref{prob}. Therefore, estimates involving norms of data and comparison between two different solutions hold true (see \cite[Corollary 5.3, Corollary 5.4 and Corollary 5.5]{SW}).

\appendix

\section{ }

In this Appendix we want to explicit the results from real analysis we use. We will begin with two well--known results which we state for the reader's convenience. The first one is a consequence of the Brezis--Lieb lemma, while the other is a generalized version of the dominated convergence theorem (see, for instance \cite{Stroock}).

\begin{Lemma}[Brezis--Lieb]
Let $\{f_n\}$ be a sequence in $L^1(\Omega)$. Then, the conditions
\begin{enumerate}
\item $f_n(x)\to f(x)$ a.e. in $\Omega$
\item $f\in L^1(\Omega)$
\item $\int_\Omega |f|=\lim_{n\to\infty}\int_\Omega |f_n|$
\end{enumerate}
imply
\[f_n\to f\quad \hbox{strongly in }L^1(\Omega)\,.\]
\end{Lemma}

\begin{Lemma}[Generalized Dominated Convergence Theorem]
Let $\{f_n \}$ and $\{g_n\}$ be sequences of measurable functions in $\Omega$. If
\begin{enumerate}
\item $g_n(x)\to g(x)$ a.e. in $\Omega$
\item $f_n\to f$ strongly in $L^1(\Omega)$
\item $|g_n|\le f_n$ for all $n\in\N$
\end{enumerate}
then
\begin{equation*}
g_n\to g\quad \hbox{strongly in }L^1(\Omega)\,.
\end{equation*}
\end{Lemma}

Given $f\in L^1(0,T)$, it follows from Lebesgue's Theorem that
\begin{equation*}
\frac1\varepsilon \int_{t-\varepsilon}^tf(s)\, ds \to f(t)\qquad \hbox{pointwise a.e. in }(0,T)\,.
\end{equation*}
Our aim in the following results is to check that we actually have strong convergence in  $L^1(0,T)$.

\begin{Lemma}\label{lema0}
Let $f\in L^1(0,T)$ be a nonnegative function with compact support. Then
\begin{equation*}
\frac1\varepsilon \int_{t-\varepsilon}^tf(s)\, ds \to f(t)\qquad \hbox{strongly in }L^1(0,T)\,.
\end{equation*}
\end{Lemma}

\begin{pf}
Define $F\colon[0,T]\to \R$ as $F(t)=\int_0^t f(s)\, ds$. This function is absolutely continuous and $F'(t)=f(t)$ a.e. Take $\varepsilon>0$ small enough to have $f(t)=0$ for all $t\in [0,\varepsilon]$. Then
\begin{multline*}
\lim_{\varepsilon\to0}\frac1\varepsilon \int_0^T\int_{t-\varepsilon}^tf(s)\, ds\, dt=
\lim_{\varepsilon\to0}\frac1\varepsilon \int_\varepsilon^T\int_{t-\varepsilon}^tf(s)\, ds\, dt
=\lim_{\varepsilon\to0}\frac1\varepsilon \int_\varepsilon^T \Big(F(t)-F(t-\varepsilon)\Big)\, dt\\
=\lim_{\varepsilon\to0}\frac1\varepsilon\left( \int_\varepsilon^TF(t)\, dt-\int_\varepsilon^T F(t-\varepsilon)\, dt \right)=\lim_{\varepsilon\to0}\frac1\varepsilon \left(\int_\varepsilon^T F(t)\, dt-\int_0^{T-\varepsilon} F(t)\, dt\right)
\\
=\lim_{\varepsilon\to0}\frac1\varepsilon \int_{T-\varepsilon}^TF(t)\, dt
=F(T)=\int_0^Tf(t)\, dt\,.
\end{multline*}
Since $\frac1\varepsilon \int_{t-\varepsilon}^tf(s)\, ds\ge 0$, it follows from the Brezis--Lieb Lemma that the convergence is in $L^1(0,T)$.
\end{pf}

\begin{Corollary}\label{corollary1}
Let $f, g\colon(0,T)\to\R$ be measurable functions such that $|g(t)|\le f(t)$ a.e. in $(0,T)$. Assume that  $f\in L^1(0,T)$  with compact support. Then
\begin{equation*}
\frac1\varepsilon \int_{t-\varepsilon}^tg(s)\, ds \to g(t)\quad \hbox{in } L^1(0,T)\,.
\end{equation*}
\end{Corollary}

\begin{pf}
Applying Lemma \ref{lema0} to $f$, we know that
\begin{equation*}
\frac1\varepsilon \int_{t-\varepsilon}^tf(s)\, ds \to f(t)\qquad \hbox{in }L^1(0,T)\,.
\end{equation*}
Our result is now  a consequence of
\begin{equation*}
\Big|\frac1\varepsilon \int_{t-\varepsilon}^t g(s)\, ds\Big|\le \frac1\varepsilon \int_{t-\varepsilon}^t f(s)\,ds
\end{equation*}
for all $\varepsilon>0$ and the generalized dominated convergence theorem.
\end{pf}

\begin{Remark}\label{RIntro}\rm
It is worth noting that similar arguments apply to functions depending on more variables. This fact allows us to justify the assertion stated in the introduction.

Let $\eta\in C_0^\infty(0,T)$ be a nonnegative function and assume $u\in L^1(0,T;  W^{1,1}(\Omega))$ and $\z\in L^\infty((0,T)\times\Omega)$. Then the function
\begin{equation*}
s\mapsto \eta(s)  \int_\Omega |\nabla u(s,x)|\, dx
\end{equation*}
belongs to $L^1(0,T; \R^N)$. By Lemma \ref{lema0},
\begin{equation*}
\lim_{\varepsilon\to0}\int_0^T \frac1\varepsilon \int_{t-\varepsilon}^t\eta(s)  \int_\Omega |\nabla u(s,x)|\, dx\, ds\, dt=\int_0^T \eta(t)  \int_\Omega |\nabla u(t,x)|\, dx\, dt\,.
\end{equation*}
Since $\displaystyle\lim_{\varepsilon\to0}\frac1\varepsilon \int_{t-\varepsilon}^t\eta(s)  |\nabla u(s,x)|\, ds=\eta(t)  |\nabla u(t,x)|$ for almost all $(t,x)\in (0,T)\times\Omega$, it follows from the nonnegativeness of all integrands that
\begin{equation}\label{ec.ap0}
\frac1\varepsilon \int_{t-\varepsilon}^t\eta(s)  |\nabla u(s,x)|\, ds\to \eta(t)  |\nabla u(t,x)|\qquad \hbox{strongly in } L^1((0,T)\times\Omega)\,.
\end{equation}
On the other hand, we also have
$\displaystyle\lim_{\varepsilon\to0}\frac1\varepsilon \int_{t-\varepsilon}^t\eta(s)  \nabla u(s,x)\, ds=\eta(t)  \nabla u(t,x)$ for almost all $(t,x)\in (0,T)\times\Omega$, owing to $\eta(s)   \nabla u(s,x)\in L^1((0,T)\times\Omega)$, wherewith
\begin{equation*}
\lim_{\varepsilon\to0}\frac1\varepsilon \int_{t-\varepsilon}^t\eta(s) \,\z(t,x)\cdot  \nabla u(s,x)\, ds=\eta(t)\, \z(t,x)\cdot  \nabla u(t,x)
\end{equation*}
for almost all $(t,x)\in (0,T)\times\Omega$. Observing that
$|\eta(s)\,\z(t,x)\cdot\nabla u(s,x)|\le \eta(s) \|\z\|_\infty|\nabla u(s,x)|$ and applying \eqref{ec.ap0}, the generalized dominated convergence theorem leads to
\begin{equation*}
\lim_{\varepsilon\to0}\int_0^T \frac1\varepsilon \int_{t-\varepsilon}^t\eta(s)  \int_\Omega \z(t,x)\cdot \nabla u(s,x)\, dx\, ds\, dt=\int_0^T \eta(t)  \int_\Omega \z(t,x)\cdot \nabla u(t,x)\, dx\, dt\,.
\end{equation*}
\end{Remark}

\begin{Lemma}\label{lema1}
Let $f\in L^1(0,T)$ with compact support and let $g\colon(0,T)\times(0,T)\to\R$ be a measurable function such that $|g(t,s)|\le f(s)$ for almost every $s,t\in(0,T)$. If
\begin{equation*}
\lim_{\varepsilon\to0}\int_0^T\frac{1}{\varepsilon}\int_{t-\varepsilon}^{t} g(t,s)\,ds\,dt = \int_0^T f(t)\,dt\,,
\end{equation*}
then
\begin{equation*}
\lim_{\varepsilon\to0}\int_0^T\frac{1}{\varepsilon}\int_{t-\varepsilon}^{t} |g(t,s)-f(t)|\,ds\, dt =  0\,,
\end{equation*}
\end{Lemma}

\begin{pf}
We define the auxiliary functions
\begin{equation*}
f_\varepsilon(t)=\dfrac{1}{\varepsilon}\int_{t-\varepsilon}^{t} f(s)\,ds\quad \hbox{ and } \quad g_\varepsilon(t)=\dfrac{1}{\varepsilon}\int_{t-\varepsilon}^{t} g(t,s)\,ds\,.
\end{equation*}
By Lemma \ref{lema0}, $f_\varepsilon\to f$ strongly in $L^1(0,T)$. Hence,
\begin{equation*}
\lim_{n\to\infty}\int_0^T \big(f_\varepsilon(t)-g_\varepsilon(t)\big)\,dt = 0\,.
\end{equation*}
Observe that almost all $t\in (0,T)$ satisfy $|g(t,s)|\le f(s)$ for almost all $s\in(0,T)$, so that
$g_\varepsilon(t)\le f_\varepsilon(t)$ holds a.e. It follows that $f_\varepsilon-g_\varepsilon\to0$ strongly in $L^1(0,T)$ and so
\begin{equation*}
g_\varepsilon\to f \qquad \hbox{strongly in }L^1(0,T)\,.
\end{equation*}
Notice also that for almost all $t\in(0,T)$, the inequality
\begin{equation*}
|g(t,s)-f(t)|\le |f(s)-f(t)|+f(s)-g(t,s)\quad\hbox{holds for almost all }s\,.
\end{equation*}
Therefore, almost all $t\in (0,T)$ satisfy
\begin{equation}\label{ec:ap1}
\frac{1}{\varepsilon}\int_{t-\varepsilon}^{t} |g(t,s)-f(t)|\,ds
\le \frac{1}{\varepsilon}\int_{t-\varepsilon}^{t} |f(s)-f(t)|\,ds+\frac{1}{\varepsilon}\int_{t-\varepsilon}^{t} f(s)-g(t,s)\,ds\,.
\end{equation}
Our aim is to check that the left hand side of \eqref{ec:ap1} tends to 0 strongly in $L^1(0,T)$. To this end, we analize its right hand side. The first term goes to 0 pointwise for almost all $t\in (0,T)$ due to Lebesgue's Theorem. Moreover, we also have the estimate
\begin{equation*}
\frac{1}{\varepsilon}\int_{t-\varepsilon}^{t} |f(s)-f(t)|\,ds\le \frac{1}{\varepsilon}\int_{t-\varepsilon}^{t} f(s)\, ds+f(t)\,.
\end{equation*}
Thus, we deduce from the generalized dominated converge theorem that
\begin{equation*}
\frac{1}{\varepsilon}\int_{t-\varepsilon}^{t} |f(s)-f(t)|\,ds\to 0\quad\hbox{strongly in }L^1(0,T)\,.
\end{equation*}
Since we already know that the second term on the right hand side of \eqref{ec:ap1} tends to 0 strongly in $L^1(0,T)$, it follows that
\[\frac{1}{\varepsilon}\int_{t-\varepsilon}^{t} |g(t,s)-f(t)|\,ds\to 0\quad\hbox{strongly in }L^1(0,T)\]
as desired.
\end{pf}

\begin{Remark}\rm
As a consequence of the previous result, we can find a subsequence
\[\frac{1}{\varepsilon_n}\int_{t-\varepsilon_n}^{t} |g(t,s)-f(t)|\,ds\]
which converges to 0 a.e. Hence,
\[\liminf_{\varepsilon\to0}\frac{1}{\varepsilon}\int_{t-\varepsilon}^{t} |g(t,s)-f(t)|\,ds=0\quad \hbox{a.e.}\]
Nevertheless, we are not able to check that every subsequence tends to 0 a.e. and so we cannot deduce that the approximate limit of $g(t.\cdot)$ at $t$ is $f(t)$.
\end{Remark}

\section*{Funding}
S. Segura de Le\'on has been supported  by MCIyU \& FEDER, under project PGC2018--094775--B--I00 and by CECE (Generalitat Valenciana) under project AICO/2021/223.

\end{document}